\documentclass[a4paper,11pt,reqno]{amsart}

\usepackage{amssymb,amsmath,array}

\makeatletter
\@addtoreset{equation}{section}

\makeatother

\def\Proof{{\bf Proof. }}
\def\endproof{\hfill$\square$}

\makeatletter
\def\follows{ \@ifnextchar[{\followsa}{\dashrightarrow}}
\def\followsa[#1]{\stackrel{#1}{\dashrightarrow}}
\makeatother

\def\smallsum{\mathop\Sigma}
\newcommand{\cone}{\mathop{\rm cone} \nolimits}
\newcommand{\conc}{\mathbin{\rm*}}

\renewcommand{\labelenumi}{{\rm \theenumi}}
\renewcommand{\theenumi}{{\rm(\roman{enumi})}}

\renewcommand{\emptyset}{\varnothing}
\renewcommand{\epsilon}{\varepsilon}
\renewcommand\phi{\varphi}

\newtheorem{theorem}{Theorem}[section]
\newtheorem{proposition}[theorem]{Proposition}
\newtheorem{lemma}[theorem]{Lemma}
\newtheorem{corollary}[theorem]{Corollary}
\newtheorem{definition}[theorem]{Definition}
\newtheorem{property}[theorem]{Property}

\renewcommand{\le}{\leqslant}
\renewcommand{\ge}{\geqslant}
\def\im{\mathop{\rm Im}}

\def\Z{\mathbf Z}
\def\Q{\mathbf Q}
\def\K{\mathbf K}
\def\U{\mathbb U}
\def\F{\mathbb F}
\def\E{\mathbb E}
\def\H{\mathbb H}
\def\X{\mathbb X}
\def\I{\mathbb I}
\def\J{\mathbb J}
\def\S{\mathbb S}
\def\C{\mathbb C}
\def\T{\mathbb T}
\def\B{\mathbb B}

\newcommand{\GL}{\mathrm{GL}}
\newcommand{\D}{\mathrm{D}}
\newcommand{\Triang}{\mathrm{T}}

\newcommand{\diag}{\mathop{\rm diag} \nolimits}
\def\c{\mathrm{cf}}
\def\soc{\mathrm{soc}}

\title[Generalization of lowering operators for ${\rm GL}_n$]
{Generalization of modular lowering operators for ${\rm GL}_n$}

\author{Vladimir Shchigolev}

\address{
\noindent
Department of Algebra, Faculty of Mathematics,
Lomonosov Moscow State University,
Leninskiye Gory, Moscow,
119899, RUSSIA}

\email{shchigolev\_vladimir@yahoo.com}


\begin{document}

\begin{abstract}
We consider the generalization of Kleshchev's lowering operators
obtained by raising all the Carter-Lusztig operators in their definition
to a power less than the characteristic of the ground field.
If we apply such an operator to a nonzero
$\GL_{n-1}$-high weight
vector
of an irreducible representation of $\GL_n$,
shall we get a nonzero
$\GL_{n-1}$-high weight
vector again?
The present paper gives the explicit answer to this question.
In this way we obtain a new algorithm for generating
some normal weights.
\end{abstract}

\maketitle

\section{Introduction}\label{Introduction}

Informally speaking, the aim of this paper is to
suggest a possible generalization of the operators
introduced by A.~S.\,Kleshchev in~\cite[Definition~2.5]{Kleshchev2}
that would be suitable for removing several nodes instead of one.
Our main sources of inspiration are therefore~\cite{Kleshchev2}
and~\cite{Carter_Lusztig}. In the latter paper, Carter and Lusztig
developed useful formulae to work with powers of their lowering
operators. We have made trivial reformulations of their results
(see Propositions~\ref{proposition:low2:1}
and~\ref{proposition:int:1}) and use them as our principal tool.

In what follows, $L_n(\lambda)$\label{Ln(lm)}
denotes the irreducible rational $\GL_n$-module
with highest weight $\lambda$.
As is well known, Kleshchev's lowering operators are enough for constructing:
\begin{itemize}
\item all $\GL_{n-1}$-high weight vectors belonging to the first level of $L_n(\lambda)$\linebreak (\cite[Theorem~4.2]{Kleshchev2});
\item all $\GL_{n-1}$-high weight vectors if $\lambda$ is a generalized Jantzen-Seitz weight\linebreak (\cite[Main Theorem]{Kleshchev_gjs11}).
\end{itemize}
The second result was proved by successive application of Kleshchev's operators to
the $\GL_{n-1}$-high weight vectors of $L_n(\lambda)$ already obtained.
However, it follows from~\cite{Shchigolev14} and the tables at the end
of~\cite{Brundan_operators} that some $\GL_{n-1}$-high weight vectors
(even belonging to levels with number less than the characteristic of the base field)
can not be reached in this way. This fact forces us to consider new lowering operators
$T_{i,j}^{(d)}(M,1)$ (see Definition~\ref{definition:low5:2}, where $R=1$),
whose effectiveness can be demonstrated by the number of $\GL_{n-1}$-high weight vectors
(up to proportionality) that can be reached by them (see Table~1).
In this paper, we exploit the simplest approach to constructing such operators:
we raise all the Carter-Lusztig operators in the definition of Kleshchev's operators
to a fixed power $d$. Thus all Kleshchev's operators correspond to $d=1$.

It would be natural to expect these new operators to behave similarly
to original Kleshchev's operators, which
they in most cases actually do.
The major obstacle
to overcome, however, is the impossibility of direct generalization
of Lemmas~2.13 and 2.15 from~\cite{Kleshchev2}, which say how the lowering
operators behave under multiplication by strictly positive elements of the hyperalgebra.
We have achieved the required modifications in Lemmas~\ref{lemma:low5:3} and~\ref{lemma:low5:4},
but at the price of introducing the extra parameter $R$ in Definition~\ref{definition:low5:2}.

We now pass to strict formulations. Let $\K$ be an algebraically closed field of
characteristic $p>0$. We assume that $\Z/p\Z\subset\K$.
We denote by $\GL_n$ the general linear group of size $n$ over $\K$.
Let $\D_n$ and $\Triang_n$ denote the subgroups of $\GL_n$ consisting
of all diagonal matrices and all upper triangular matrices respectively.
We put $X(n):=\Z^n$\label{X(n)}
and call the elements of this set {\it weights}.
Any weight $(\lambda_1,\ldots,\lambda_n)$ will be identified
with the character of $\D_n$ that takes $\diag(t_1,\ldots,t_n)$ to
$t_1^{\lambda_1}\cdots t_n^{\lambda_n}$. A vector $v$ of a rational $\GL_n$-module is called
a {\it $\GL_n$-high weight vector} if the line $\K\cdot v$
is fixed by $\Triang_n$. We denote by $X^+(n)$ the subset of $X(n)$ consisting of
all weakly decreasing sequences and call the elements of $X^+(n)$\label{domX(n)}
{\it dominant weights}.

For $n>1$, the group $\GL_{n-1}$ will be identified with a subgroup of $\GL_n$
consisting of matrices having zeros in the last row and the last column
except the position of their intersection, where they have $1$.

The main results of this paper are as follows. For a $\GL_{n-1}$-high weight vector
$v$ of the irreducible module $L_n(\lambda)$,
{
\renewcommand{\theenumi}{{\rm(\arabic{enumi})}}
\begin{enumerate}
    \item\label{main_result:1}
          Theorems~\ref{theorem:commpoly:1} and~\ref{theorem:commpoly:2} show if
          the vector $T_{i,j}^{(d)}(M,1)v$ is a nonzero
          $\GL_{n-1}$-high weight vector for any $1\le d<p$, $1\le i<j\le n$
          and $M\subset(i..j)$;
     \item\label{main_result:2}
           Theorems~\ref{theorem:crit:0.5} and~\ref{theorem:crit:1}\ref{theorem:crit:1:part:2}
           show if
           for fixed $1\le d<p$ and $1\le i<j\le n$, there exists $M$
           such that $T_{i,j}^{(d)}(M,1)v$ is a nonzero
           $\GL_{n-1}$-high weight vector.
\end{enumerate}
}

Of main interest for us, however, are the weights $\mu\in X^+(n-1)$ such that
there exists a nonzero $\GL_{n-1}$-high weight vector of $L_n(\lambda)$ having weight $\mu$.
These weights are called {\it normal}. It is a major problem to find
a direct combinatorial description of all normal weights,
from the solution of which the structure of the socle of
the restriction $L_n(\lambda)\downarrow_{\GL_{n-1}}$ would follow
(see~\cite[Theorem D]{Kleshchev_bou}).
Now we are going to formulate a property of the normal weight pattern
following from result~\ref{main_result:2} above.
This property implies an algorithm to construct normal weights
that generalizes any such algorithm
of~\cite{Kleshchev2},~\cite{Kleshchev_gjs11} or~\cite{Shchigolev14}.

\begin{definition}\label{def:introduction:1}
Let $A\subset\Z^n$, $B\subset\Z^m$ and
$\varphi:A\to B$. We say that $\varphi$ is weakly increasing
{\rm(}decreasing{\rm)}
w.r.t. the $k$th coordinate, where $1\le k\le\min\{n,m\}$ if
$y_k{\ge}x_k$
{\rm(}resp. $y_k{\le}x_k${\rm)}, whenever
$\phi(x_1,\ldots,x_n)=(y_1,\ldots,y_m)$.
\end{definition}
For integers $1\le i<j\le n$, $1\le d<p$ and weights $\mu\in X^+(n-1)$, $\lambda\in X^+(n)$,
we put (assuming $j<n$ in the second line)
{\arraycolsep=2pt
$$\label{CXX}
\begin{array}{rcl}
\mathfrak C^\mu(i,j)&:=&{\{}t\in(i..j):t-i+\mu_i-\mu_t\equiv0\!\!\!\!\pmod p\},\\[6pt]
\mathfrak X^\mu_d(i,j)    &:=&{\{} (t,s)\in [i..j)\times[1..d]: t-i+\mu_i-\mu_{t+1}\equiv d-s\!\!\!\!\pmod p{\}},\\[6pt]
\mathfrak X^{\mu,\lambda}_d(i,j)&:=&{\{} (t,s)\in [i..j)\times[1..d]:t-i+\mu_i-\lambda_{t+1}\equiv d-s\!\!\!\!\pmod p {\}}.
\end{array}
$$
}
In these definitions
and in what follows, by $[i..j]$, $[i..j)$, $(i..j]$, $(i..j)$\label{intervals}
we denote the sets $\{x\in\Z:i\le x\le j\}$, $\{x\in\Z:i\le x<j\}$,$\{x\in\Z:i<x\le j\}$,
$\{x\in\Z:i<x<j\}$ respectively. Moreover, for any sequence $a$, we denote
its length by $|a|$ and its $i$th entry by $a_i$.

\begin{property}
Let $\mu\in X^+(n-1)$ be normal for $\lambda\in X^+(n)$ and $1\le d<p$.
\begin{enumerate}
    \item If $1\le i<n$ and there exists an
    injection
    $\epsilon:\mathfrak X^{\mu,\lambda}_d(i,n)\to\mathfrak C^\mu(i,n)$
    weakly decreasing w.r.t the first coordinate, then the weight
    $(\mu_1,\ldots,$\linebreak $\mu_{i-1},\mu_i-d,\mu_{i+1},\ldots,\mu_{n-1})$
    is also normal for $\lambda$.
    \item If $1\le i<j<n$, $(j-1,1)\in\mathfrak X^\mu_d(i,j)$,
             $(j-1,1)\notin\mathfrak X^{\mu,\lambda}_d(i,j)$,
    there exists an
    injection
    $\epsilon:\mathfrak X^{\mu,\lambda}_d(i,j)\to\mathfrak C^\mu(i,j)$
    weakly decreasing w.r.t the first coordinate and an
    injection
    $\tau:\mathfrak X^{\mu,\lambda}_d(i,j)\to \mathfrak X^\mu_d(i,j)\setminus\{(j-1,1)\}$
    weakly increasing w.r.t. the first coordinate and
    weakly decreasing w.r.t. the second coordinate, then
    the weight $(\mu_1,\ldots,\mu_{i-1},\mbox{$\mu_i-d$},\mu_{i+1},\ldots,\mu_{j-1},$ $\mu_j+d,\mu_{j+1},\ldots,\mu_{n-1})$
    is also normal for $\lambda$.
\end{enumerate}
\end{property}

The table below shows the maximum over all $p$-restricted weights
$\lambda\in X^+(n)$ of
the number of the $\GL_{n-1}$-high weight vectors
in $L_n(\lambda)$ (up to proportionality) that can be reached
by the operators $T_{i,j}^{(d)}(M,1)$ for $1\le d<p$,
but can not be reached by Kleshchev's lowering operators
(i.e. $T_{i,j}^{(1)}(M,1)$).

\medskip
\begin{center}
{\small
\begin{tabular}{|c||c|c|c|c|c|c|c||c|c|c|c|c||c|c|c|c||}
\hline
$p$&\multicolumn{7}{c||}{3}&\multicolumn{5}{c||}{5}&\multicolumn{4}{c||}{7}\\
\hline
$n$&2&3&4& 5& 6&  7&  8&2 &3& 4&  5&6  &2&3& 4&  5\\
\hline
max&0&1&3&14&43&153&465&0 &4&26&153&917&0&9&87&713\\
\hline
\end{tabular}\\[4pt]
Table 1
}
\end{center}

\medskip

We fix the following sequences of $X(n)$: $\epsilon_i$ having $1$
at position $i$, where $1\le i\le n$, and zeros elsewhere;
$\alpha(i,j):=\epsilon_i-\epsilon_j$, where $1\le i,j\le n$;
$\alpha_i=\alpha(i,i+1)$, where $1\le i<n$.\label{weights}
The lengths of these sequences will always be clear from the context
and the following stipulation: an elementwise linear combination of
two sequences $a$ and $b$ is well defined only if $|a|=|b|$.
Elements of $X(n)$ are ordered as follows: $\lambda\ge\mu$ if
$\lambda-\mu$ is a sum of $\alpha_i$ with nonnegative coefficients.
The descending factorial power $x^{\underline n}$\label{descending_factorial}
equals $x(x-1)\cdots(x-n+1)$
if $n\ge0$ and equals $1/(x+1)\cdots(x-n)$ if $n<0$. The principal relation
for this power is $x^{\underline{m+n}}=x^{\underline m}(x-m)^{\underline m}$.
Any product of the form $\prod_{i=a}^bX_i$, with not necessarily commuting
factors $X_i$, will mean $X_a\cdots X_b$. Following the standard
agreement, we shall interpret any expression $a^n$ as the sequence of
length $n$ whose every entry is $a$, if this notation does not cause
confusion. For example, $(a^3,1^2)=(a,a,a,1,1)$. For two finite sequences
$a=(a_1,\ldots,a_n)$ and $b=(b_1,\ldots,b_m)$, we define their product by
$a*b:=(a_1,\ldots,a_n,b_1,\ldots,b_m)$.
A formula $A\sqcup B=C$ will mean $A\cup B=C$ and $A\cap B=\emptyset$.
For any condition $\mathcal P$, let $\delta_{\mathcal P}$ be $1$
if $\mathcal P$\label{delta}
is true and $0$ if it is false.

The layout of the paper is as follows. Section~\ref{Hyperalgebra} is devoted
to various basic constructions in the hyperalgebra over integers $\U(n)$
and its extension $\bar\U(n)$. The central
topics here are Lemma~\ref{lemma:algebras:2} on intersections and
Lemmas~\ref{lemma:low5:1} and~\ref{lemma:low5:2} on block products
(see also Definition~\ref{definition:low5:1}).
In Section~\ref{Multiplication formulae}, we introduce the operators
$T_{i,j}^{(d)}(M,R)$ (Definition~\ref{definition:low5:2}) and
prove the multiplication formulae for them
(Lemmas~\ref{lemma:low5:3} and~\ref{lemma:low5:4}).
In Section~\ref{Coefficients}, we introduce the rational expressions
similar to $\xi_{r,s}(M)$ of~\cite[Definition 2.9]{Kleshchev2}
virtually for the same purpose (see~(\ref{equation:coeff:2.5}),
(\ref{equation:coeff:2.75}) and Lemma~\ref{lemma:coeff:1}).
While up to now we have not been interested whether
our elements belong to $\U(n)$, Section~\ref{Integral elements}
considers this question. Finally, in Section~\ref{Proof of the main results},
the main results are proved.

The author would like to thank J.\;Brundan, who kindly
made his GAP program (calculating normal and good weights)
available to him.

\section{Hyperalgebra over $\Z$}\label{Hyperalgebra}

In Sections~\ref{Hyperalgebra}--\ref{Integral elements},
we fix an integer $n>0$.
Let $U_\Q(n)$\label{universal_enveloping_algebra}
be the universal enveloping algebra of
the Lie algebra ${\mathfrak gl}_\Q(n)$ of all $n\times n$-matrices
over rationals.
As usual, ${\mathfrak gl}_\Q(n)$ is embedded into $U_\Q(n)$.
We denote by $\X_{i,j}$ the matrix of ${\mathfrak gl}_\Q(n)$ having $1$
at $ij$-entry and $0$ elsewhere.
The ring $U_\Q(n)$ is graded by the subgroup of $\Z^n$ consisting of sequences
with sum zero so that $\X_{i,j}$ has weight $\alpha(i,j)$.
In this paper,
weights and homogeneity of elements are always understood
with respect to this grading.

The {\it hyperalgebra} over $\Z$ is the subring $\U(n)$\label{hyperalgebra}
of $U_\Q(n)$ generated by
$$
\begin{array}{rcll}
\X_{i,j}^{(r)}&:=&\frac{(\X_{i,j})^r}{r!}&\mbox{ for integers }1\le i\ne j\le n \mbox{ and }r\ge0;\\[12pt]
\binom{\X_{i,i}}{r}&:=&\frac{\X_{i,i}(\X_{i,i}-1)\cdots(\X_{i,i}-r+1)}{r!}&\mbox{ for integers }1\le i\le n\mbox{ and }r\ge0.
\end{array}
$$
In this definition, empty product means the identity of $U_\Q(n)$.
It is convenient to define the above elements
as zero
if $r$ is a negative integer. For mutually distinct $i,j,k$, we have
{
\begin{equation}\label{equation:algebras:0.5}
\arraycolsep=0pt
\begin{array}{l}
\displaystyle \X_{i,j}^{(a)}\X_{j,i}^{(b)}=\sum_t\X_{j,i}^{(b-t)}\binom{\H_i-\H_j-(a+b)+2t}{t}\X_{i,j}^{(a-t)},\\[12pt]
\displaystyle \X_{i,j}^{(a)}\X_{j,k}^{(b)}{=}\sum_t\X_{j,k}^{(b-t)}\X_{i,k}^{(t)}\X_{i,j}^{(a-t)},\;\;
\displaystyle \X_{j,k}^{(b)}\X_{i,j}^{(a)}{=}\sum_t(-1)^t\X_{i,j}^{(a-t)}\X_{i,k}^{(t)}\X_{j,k}^{(b-t)}.\\[12pt]
\end{array}
\end{equation}
}

We shall use the notation $\E^{(r)}_{i,j}:=\X^{(r)}_{i,j}$,
$\F^{(r)}_{i,j}:=\X^{(r)}_{j,i}$, where $1\le i<j\le n$,
and $\H_i:=\X_{i,i}$.
We also put
$\E^{(r)}_i:=\E^{(r)}_{i,i+1}$ and omit the superscript ${}^{(1)}$.

Let $UT(n)$\label{UT(n)}
denote the set of integer $n\times n$-matrices $N$
with nonnegative entries such that $N_{a,b}=0$ unless $a<b$ and $e_{i,j}$,
where $i<j$, denote the element of $UT(n)$ with $1$ at $ij$-entry and $0$
elsewhere.
For any integer $n\times n$-matrix $N$, we define
$$\label{FNEN}
\F^{(N)}:=\prod_{1\le a<b\le n}\F_{a,b}^{(N_{a,b})},\quad
\E^{(N)}:=\prod_{1\le a<b\le n} \E_{a,b}^{(N_{a,b})},
$$
where $\F_{a,b}^{(N_{a,b})}$ precedes $\F_{c,d}^{(N_{c,d})}$
iff $b<d$ or $b=d$ and $a<c$ in the first product and
$\E_{a,b}^{(N_{a,b})}$ precedes $\E_{c,d}^{(N_{c,d})}$
iff $a<c$ or $a=c$ and $b<d$ in the second product.
Obviously, $\F^{(N)}=\E^{(N)}=0$ if $N$ contains a negative entry.

\begin{proposition}[\mbox{\cite[2.1]{Carter_Lusztig}}]\label{proposition:algebras:1}
Elements
$
\textstyle
\F^{(N)}\binom{\H_1}{r_1}\cdots\binom{\H_n}{r_n}\E^{(M)},
$
where $N,M$\linebreak$\in UT(n)$ and $r_1,\ldots,r_n$ are nonnegative integers
form a $\Q$-basis of $U_\Q(n)$.
These elements generate $\U(n)$ as a $\Z$-module.
\end{proposition}
In particular, the $\Q$-subalgebra $U^0_\Q(n)$\label{U0}
of $U_\Q(n)$
generated by $\H_1,\ldots,\H_n$ is generated by them
freely as a commutative $\Q$-algebra.

For any $N\in UT(n)$, we define the ring automorphism $\tau_N$\label{tauN}
of $U^0_\Q(n)$
by $\tau_N(\H_i):=\H_i+\sum_{1\le a<i}N_{a,i}-\sum_{i<b\le n}N_{i,b}$
for $i=1,\ldots,n$. This definition is made to ensure
\begin{equation}\label{equation:algebras:1}
f\F^{(N)}=\F^{(N)}\tau_N(f),\quad\E^{(N)}f=\tau_N(f)\E^{(N)}
\end{equation}
for every $f\in U^0_\Q(n)$ and $N\in UT(n)$.
It follows from these formulae and Proposition~\ref{proposition:algebras:1}
that every element $f\in U^0_\Q(n)\setminus\{0\}$
is not a zero divisor in $U_\Q(n)$.
Let $\bar\U(n)$\label{barU(n)}
be the right ring of quotients of $U_\Q(n)$
with respect to $U^0_\Q(n)\setminus\{0\}$. Its existence follows
from the right Ore condition, which can be easily checked.
It is easy to prove that $\bar\U(n)$ is also the left ring of quotients
for the same pair. In Sections~\ref{Hyperalgebra}--\ref{Integral elements},
we use miscellaneous rings, all of which are subrings of $\bar\U(n)$.
Below we give the table that defines them.

\smallskip

{\arraycolsep=-1pt\label{miscellaneous_rings}

\setlength{\extrarowheight}{4pt}

\begin{tabular}{|l||c|c|c|}
\hline
$R$         & $\U$                              & $\bar\U$              & $\mathcal U$ \\
\hline
$R^+(a,b)$  & $\begin{array}{c} \bigl\{\E_{i,j}^{(r)}:a{\le}i{<}j{\le}b\bigr\}  \\\mbox{subring}\end{array}$  &                       &                        \\
\hline
$R^0(a,b)$  &
$\begin{array}{c}\bigl\{\binom{\H_k}r:a{\le}k{\le}b\bigl\}\\\mbox{subring}\end{array}$  &
$\begin{array}{c}\bigl\{\H_k:a{\le}k{\le}b\bigr\}\\\mbox{subfield}\end{array}$&
$\begin{array}{c}\bigl\{\H_k:a{\le}k{\le}b\bigr\}\\\mbox{subring}\end{array}$\\
\hline
$R^-(a,b)$  & $\begin{array}{c}\bigl\{\F_{i,j}^{(r)}:a{\le}i{<}j{\le}b\bigr\}\\\mbox{subring}\end{array}$  &                       &                        \\
\hline
$R^{-,0}(a,b)$  &
$\begin{array}{c}
\bigl\{\F_{i,j}^{(r)}:a{\le}i{<}j\le b\bigr\}\\
\cup\bigl\{\binom{\H_k}r:a{\le}k{\le}b\bigl\}\\
\mbox{subring}
\end{array}$  &
$\begin{array}{c}
\bigl\{\F_{i,j}^{(r)}:a{\le}i{<}j{\le}b\bigr\}\\
\cup\;\bar\U^0(a,b)\\
\mbox{subring}
\end{array}$
&
$\begin{array}{c}
\bigl\{\F_{i,j}^{(r)}:a{\le}i{<}j{\le}b\bigr\}\\
\cup\bigl\{\H_k:a{\le}k{\le}b\bigl\}\\
\mbox{subring}
\end{array}$
\\
\hline
\end{tabular}}

\smallskip

\noindent
The last three columns correspond to values of $R$,
the last four cells in each of these columns contain
generating sets and the word showing if we generate
a subring or a subfield. We abbreviate $R^?(1,n)$ to $R^?(n)$.

If we extend the automorphism $\tau_N$ of $U^0_\Q(n)$
to the automorphism of the field $\bar\U^0(n)$ by
$\tau_N(fg^{-1})=\tau_N(f)\tau_N(g)^{-1}$,
where $f,g\in U^0_\Q(n)$ and $g\ne0$,
then formulae~(\ref{equation:algebras:1})
will hold for any $f\in\bar\U^0(n)$.

\begin{lemma}\label{lemma:algebras:1}
Any element $x\in\bar\U(n)$ is uniquely represented in the form
$
x=\sum_{N,M\in UT(n)}\F^{(N)}H_{N,M}\E^{(M)},
$
where $H_{N,M}\in\bar\U^0(n)$.
\end{lemma}
\Proof It suffices to apply~(\ref{equation:algebras:1})
and Proposition~\ref{proposition:algebras:1}. \endproof

We shall use two types of ideals in either ring $\U(n)$ and $\bar\U(n)$.
For $S\subset[1..n)$, let $\I_S$ and $\bar\I_S$\label{idealsI} denote the
left ideals of $\U(n)$ and $\bar\U(n)$ respectively generated by
$\E^{(r)}_i$, where $r\ge1$ and $i\in S$.
Clearly $\bar\I_S$ is the left ideal of $\bar\U(n)$
generated by $\E_i$, where $i\in S$.
Let $C=(c_1,\ldots,c_{n-1})$ be the sequence with entries
belonging to $\Z\cup\{+\infty\}$.
Let $\J^{(C)}$ and $\bar\J^{(C)}$\label{idealsJ}
denote the left ideals of
$\U(n)$ and $\bar\U(n)$ respectively generated by
homogeneous elements having weight with $\alpha_i$-coefficient
strictly greater than $c_i$ for some $i=1,\ldots,n-1$.
Obviously
$\I_\emptyset=\bar\I_\emptyset
=\J^{((+\infty,\ldots,+\infty))}=\bar\J^{((+\infty,\ldots,+\infty))}=0$.

\begin{lemma}\label{lemma:algebras:2}
$(\bar\I_S+\bar\J^{(C)})\cap\U(n)=\I_S+\J^{(C)}$.
\end{lemma}
\Proof
Applying the automorphism of $\bar\U(n)$
permuting the indexes, it suffices to consider the case $S=[k..n)$.
We only need to prove that the left-hand side is contained in right-hand side.

Let $f_1\in\bar\I_S$.
The definition of $\bar\I_S$ and Lemma~\ref{lemma:algebras:1} yield
$$
f_1=\sum_{i=k}^{n-1}\left(\sum_{N,M\in UT(n)}\F^{(N)}H^{(i)}_{N,M}\E^{(M)}\right)\E_i,
$$
where $H^{(i)}_{N,M}\in\bar\U^0(n)$.
The order of factors in $\E^{(M)}$, we agreed upon,
and formulae~(\ref{equation:algebras:0.5})
show that $\E^{(M)}\E_i$ is equal to the $\Z$-linear combination
of elements $\E^{(M')}$, where $M'\in UT(n)$ such that
$M'_{s,t}\ne0$ for some $s$ and $t$ with $i\le s<t\le n$.

Now let $f_2\in\bar\J^{(C)}$. It is easy to see that
$
f_2=\sum_{N,M\in UT(n)}\F^{(N)}{\mathcal H}_{N,M}\E^{(M)},
$
where ${\mathcal H_{N,M}}\in\bar\U^{(0)}(n)$ and
if ${\mathcal H}_{N,M}\ne0$ then there is some $i=1,\ldots,n-1$
such that the $\alpha_i$-coefficient of the weight of
$\E^{(M)}$ is strictly greater than $c_i$.

Now assume $f_1+f_2\in\U(n)$. The above representations of
$f_1$ and $f_2$ and Lemma~\ref{lemma:algebras:1} yield
$
f_1+f_2=\sum_{N,M\in UT(n)}\F^{(N)}H_{N,M}\E^{(M)},
$
where $H_{N,M}\in\U^0(n)$ and if $H_{N,M}\ne0$ then
either the $\alpha_i$-coefficient of the weight of $\E^{(M)}$
is greater than $c_i$ for some $i=1,\ldots,n-1$ or
$M_{s,t}\ne0$ for some $s$ and $t$ with $k\le s<t\le n$.
In the first case $\E^{(M)}\in\J^{(C)}$.

Therefore, consider the second case. Let $(i,j)$ be the lexicographically
greatest pair such that $1\le i<j\le n$ and $M_{i,j}>0$.
Hence $(i,j)$ is greater than or equal to $(s,t)$ and thus $k\le s\le i<j\le n$.
Hence by~(\ref{equation:algebras:0.5}) (inductively)
we get $\E_{i,j}^{(M_{i,j})}\in\I_{[k..n)}$.
Now recall that $\E^{(M)}=\E^{(M-M_{i,j}e_{i,j})}\E_{i,j}^{(M_{i,j})}$ and
thus $\E^{(M)}\in\I_{[k..n)}$. \endproof

In what follows, we shall use the abbreviations
$$\label{CB}
\C(i,j):=j-i+\H_i-\H_j,\quad\B(i,j):=j-i+\H_i-\H_{j+1}.
$$
The main tool of our investigation is the following elements of
the hyperalgebra
$$\label{S}
\S_{i,j}:=\sum_{A\subset(i..j)}
\left(\F_{i,j}^A\prod_{t\in(i..j)\setminus A}\C(i,t)\right),
$$
where $1\le i\le j\le n$ and
$\F_{i,j}^A=\F_{a_0,a_1}\F_{a_1,a_2}\cdots\F_{a_m,a_{m+1}}$ for
$A\cup\{i,j\}=\{a_0<\cdots<a_{m+1}\}$.
In particular, $\S_{i,i}=1$.
Elements $S_{i,j}$ were first introduced in~\cite{Carter_Lusztig}.
In this connection, we call them the {\it Carter--Lusztig lowering operators}.

The operators $\S_{i,j}$ possess the property
\begin{equation}\label{equation:low2:0.5}
\E_{l-1}\S_{i,j}\equiv 0\pmod{\U(n)\cdot\E_{l-1}},
\end{equation}
for $l\ne j$ (see~\cite[Lemma~2.4]{Carter_Lusztig}
or~\cite[Lemma~3.6]{Brundan_operators}).
If $j=l$ then it follows from~\cite[Lemma 2.6]{Carter_Lusztig}
that the following equivalences hold.

\begin{proposition}\label{proposition:low2:1}
Let $1\le i<j\le n$ and $l,d\ge 0$. Then modulo the left ideal of $\U(n)$
generated by $\E_{j-1},\ldots,\E_{j-1}^{(l)}$, we have
$$
\begin{array}{l}
\E_{j-1}^{(l)}\S_{i,j}^d\equiv\binom dl\S_{i,j-1}^l\S_{i,j}^{d-l}
(\C(i,j)-d+l-1)\cdots(\C(i,j)-d)
\mbox{ if }l\le d;\\[8pt]
\E_{j-1}^{(l)}\S_{i,j}^d\equiv0\mbox{ if }l>d.
\end{array}
$$
\end{proposition}

\begin{corollary}\label{corollary:low3:1}
Let $1\le i<j_1\le\cdots\le j_d\le n$ and $l=1,\ldots,d$.
Then modulo the left ideal of $\U(n)$ generated by $\E_{j_l-1}$, we have
$$
\E_{j_l-1}\S_{i,j_1}\cdots\S_{i,j_d}\equiv
\S_{i,j'_1}\cdots\S_{i,j'_d}(d_{j_l}-d_{j_l+1})(\C(i,j_l)-d_{j_l}),
$$
where
$d_s:=|\{t\in[1..d]:s\le j_t\}|$,
$j'_t=j_t-1$ if $t$ equals the smallest of the numbers of $[1..d]$ such that $j_t=j_l$,
and $j'_t=j_t$ otherwise.
\end{corollary}

We shall abbreviate $\S_{i,J}:=\S_{i,j_1}\cdots\S_{i,j_d}$, where
$J=(j_1,\ldots,j_d)$ is a sequence of integers.
For $N\in UT(n)$ and $1<m\le n$, one can easily check the formula
{\arraycolsep=-1pt
\begin{equation}\label{equation:low3:1}
\begin{array}{l}
\displaystyle [\E_{m-1},\F^{(N)}]=\sum_{1\le s<m-1}(N_{s,m-1}+1)\F^{(N-e_{s,m}+e_{s,m-1})}\\[6pt]
\displaystyle +\F^{(N-e_{m-1,m})}\left(\H_{m-1}-\H_m+1-\sum_{m-1<b\le n}N_{m-1,b}+\sum_{m<b\le n}N_{m,b}\right)\\[6pt]
\displaystyle -\sum_{m<t\le n}(N_{m,t}+1)\F^{(N-e_{m-1,t}+e_{m,t})}.
\end{array}
\end{equation}}

For $i=1,\ldots,n-1$, let $\theta_i$\label{theta}
be the abelian group automorphism of
$\bar\U(n)$ defined by the rules:
\begin{itemize}
\itemsep=3pt
\item
$\theta_i(\H_j)=\H_j+\delta_{j=i+1}-\delta_{j=i}$ for any $j\in[1..n]$
and the restriction of $\theta_i$ to $\bar\U^0(n)$ is a field automorphism;
\item $\theta_i(FHE)=F\theta_i(H)E$, where $F\in\U^-(n)$, $H\in\bar\U^0(n)$
and $E\in\U^+(n)$.
\end{itemize}

\begin{proposition}\label{proposition:low5:-1&low5:2.5}
${}$
\begin{enumerate}
\item\label{proposition:low5:-1}
Let $1\le a\le b\le n$, $1<m\le n$ and $N\in UT(n)$ be such that
$\F^{(N)}\in\U^-(a,b)$.
Then $[\E_{m-1},\F^{(N)}]\in\mathcal U^{-,0}(a,b)$.
\item\label{proposition:low5:2.5}
Let $R\in\bar\U^{-,0}(a,b)$. We put $R_0:=R$. Suppose that every $R_{k+1}$
is obtained from $R_k$ either by left multiplication by $\E_i$,
where $i\in M$, or by application of some $\theta_j$.
Then all $R_m$ belong to $\bar\U^{-,0}(a,b)\oplus\bar\I_M$.
\end{enumerate}
\end{proposition}
\Proof
 This result follows immediately from~(\ref{equation:low3:1}).
\endproof

For $1 \le a\le b\le c\le d\le n$, we have
\begin{equation}\label{equation:low5:0.5}
[\U^\pm(a,b),\U^\mp(c,d)]=0.
\end{equation}

\begin{proposition}\label{proposition:low5:0}${}$
\begin{enumerate}
\itemsep=3pt
\item\label{proposition:low5:0:part:1}
      If $1\le a\le b\le i<n$ and $Y\in\U^-(a,b)\bar\U^0(n)$, then
      $\E_iY=\theta_i(Y)\E_i$.
\item\label{proposition:low5:0:part:2}
      If $1\le i<n$, $m\in\Z$,
      $X\in\bar\U(n)$ and $H\in\bar\U^0(n)$, then
      $\theta_i^m(XH)=\theta_i^m(X)\theta_i^m(H)$.
\item\label{proposition:low5:0:part:3}
      If $1{\le}i{<}n$, $m{\in}\Z$ and $Y_1,Y_2\in\bar\U^{-,0}(n)$,
      then $\theta_i^m(Y_1Y_2)=\theta_i^m(Y_1)\theta_i^m(Y_2)$.
\item\label{proposition:low5:0:part:4}
      If $1\le i<j<n$ and $m\in\Z$, then
      $\theta_j^m\bigl(\C(i,j)\bigl)=\C(i,j)+m$.
\item\label{proposition:low5:0:part:5}
      The ideals $\I_M$ and $\bar\I_M$ are stable under $\theta_i$.
\end{enumerate}
\end{proposition}
\Proof It suffices to notice that the restriction of $\theta_i$ to
$\bar\U^0(n)$ commutes with $\tau_N$ and apply
formulae~(\ref{equation:algebras:1}) and~(\ref{equation:low5:0.5}).
\endproof

\begin{lemma}\label{lemma:low5:1}
Suppose $1\le a_1<b_1\le\cdots\le a_k<b_k\le n$
and for each $i=1,\ldots,k$, we are given an element
$X_i{\in}\bar\U^{-,0}(a_i,b_i)\oplus\bar\I_{M_i}$,
where $M_i\subset[a_i..b_i)$.
\begin{enumerate}
\item\label{lemma:low5:1:part:1} If for each $i=1,\ldots,k$,
     we are given another $X'_i\in\bar\U(n)$ such that
     $X_i\equiv X'_i\pmod{\bar\I_{M_i}}$, then
     $X_1\cdots X_k\equiv X'_1\cdots X'_k\pmod{\bar\I_{M_1\cup\cdots\cup M_k}}$.
\item\label{lemma:low5:1:part:2} If $s\in[1..n)$ then modulo
      $\bar\I_{M_1\cup\cdots\cup M_k\cup\{s\}}$ the element
      $\E_sX_1\cdots X_k$ equals
      {
      \renewcommand{\theenumi}{}
      \begin{enumerate}
      \item\label{lemma:low5:1:case:1} $0$ if $s\notin[a_1..b_1)\cup\cdots\cup[a_k..b_k)${\rm ;}
      \item\label{lemma:low5:1:case:2} $X_1\cdots X_{m-1}\E_sX_m\cdots X_k$ if $s\in[a_m..b_m)$ and $m=1$ or $s>b_{m-1}${\rm ;}
      \item\label{lemma:low5:1:case:3} $X_1\cdots X_{m-2}\theta_s(X_{m-1})\E_sX_m\cdots X_k$ if
      $m>1$ and $s=a_m=b_{m-1}$.
      \end{enumerate}
      }
\end{enumerate}
\end{lemma}
\Proof
Let $Y_i$ be the element of $\bar\U^{-,0}(a_i,b_i)$ such that
$X_i\equiv Y_i\pmod{\bar\I_{M_i}}$.

\ref{lemma:low5:1:part:1} We use induction on $k$. The case $k=0$
is obvious. Assume that $k>0$ and that the assertion is true for
smaller number of factors.
By the inductive hypothesis, we have
$X_2\cdots X_k\equiv X'_2\cdots X'_k\pmod{\bar\I_{M_2\cup\cdots\cup M_k}}$,
whence $X'_1X_2\cdots X_k\equiv X'_1X'_2\cdots X'_k\pmod{\bar\I_{M_2\cup\cdots\cup M_k}}$.
Therefore, it remains to prove that
$(X_1-X'_1)X_2\cdots X_k\in\bar\I_{M_1\cup\cdots\cup M_k}$.
Applying the inductive hypothesis one more time, we get
$X_2\cdots X_k\equiv Y_2\cdots Y_k\pmod{\bar\I_{M_2\cup\cdots\cup M_k}}$.
Hence $(X_1-X'_1)X_2\cdots X_k\equiv (X_1-X'_1)Y_2\cdots Y_k\pmod{\bar\I_{M_2\cup\cdots\cup M_k}}$.
Moreover, (\ref{equation:low5:0.5}) implies
$(X_1-X'_1)Y_2\cdots Y_k\in \bar\I_{M_1}$.

\ref{lemma:low5:1:part:2}
Part~\ref{lemma:low5:1:part:1} yields
$\E_sX_1\cdots X_k\equiv \E_sY_1\cdots Y_k\pmod{\bar\I_{M_1\cup\cdots\cup M_k}}$.
In\linebreak case~\ref{lemma:low5:1:case:1}, the result immediately
follows from~(\ref{equation:low5:0.5}).

Now we consider the case where $s\in[a_m..b_m)$ for some $m=1,\ldots,k$.
Then Proposition~\ref{proposition:low5:0}\ref{proposition:low5:0:part:1}
implies $\E_sY_1\cdots Y_k=\theta_s(Y_1)\cdots\theta_s(Y_{m-1})\E_sY_m\cdots Y_k$.
Note that $\E_sY_m\cdots Y_k\in\bar\U^{-,0}(a_m,b_k)\oplus\bar\I_{\{s\}}$,
which easily follows from Proposition~\ref{proposition:low5:-1&low5:2.5}\ref{proposition:low5:-1}.
Let us prove inductively on $j=0,\ldots,m-1$ that
$$
\begin{array}{l}
\theta_s(Y_1)\cdots\theta_s(Y_{m-1})\E_sY_m\cdots Y_k\equiv
\theta_s(X_1)\cdots\theta_s(X_j)\times\\[6pt]
\times\theta_s(Y_{j+1})\cdots\theta_s(Y_{m-1})\E_sY_m\cdots Y_k
\pmod{\bar\I_{M_1\cup\cdots\cup M_j\cup\{s\}}}.
\end{array}
$$
Actually the proof follows from
$\theta_s(Y_{j+1})\equiv \theta_s(X_{j+1})\pmod{\bar\I_{M_{j+1}}}$,
which holds by
Proposition~\ref{proposition:low5:0}\ref{proposition:low5:0:part:5}.
The equivalence we just proved applied for $j=m-1$ and
part~\ref{lemma:low5:1:part:1} of the current lemma yield
$$
\begin{array}{l}
\E_sX_1\cdots X_k\equiv
\theta_s(X_1)\cdots\theta_s(X_{m-1})\E_sY_m\cdots Y_k\equiv \\[6pt]
\theta_s(X_1)\cdots\theta_s(X_{m-1})\E_sX_m\cdots X_k
\pmod{\bar\I_{M_1\cup\cdots\cup M_k\cup\{s\}}}.
\end{array}
$$
Cases~\ref{lemma:low5:1:case:2} and~\ref{lemma:low5:1:case:3}
now follow from the definition of $\theta_s$.
\endproof

We shall use the notation $\E(i,j):=\E_i\cdots\E_{j-1}$, where
$1\le i\le j\le n$.

\begin{definition}\label{definition:low5:1}
Let $1\le j\le n$ and $K=(k_1,\ldots,k_m)$, $L=(l_1,\ldots,l_q)$
be sequences of integers such that $1\le k_1\le\cdots\le k_m\le j$ and
$q=m$ or $q=m+1$. For $f\in\bar\U(n)$, we define $\E_j(K,L)(f)$ to equal
\begin{itemize}
\item $f$ if $K=L=\emptyset$;
\item $\E(k_m,j)\cdot\E_j(K',L)(f)$
  if $m=q>0$, where $K'=(k_1,\ldots,k_{m-1})$;
\item $\theta_j^{l_{m+1}}\bigl(\E_j(K,L')(f)\bigl)$
  if $q=m+1$, where $L'=(l_1,\ldots,l_m)$.
\end{itemize}
\end{definition}
The maps $\E_j(K,L)$ we have defined are
abelian group automorphisms of $\bar\U(n)$.
We also see that $\E_j(K,L)=\E_j(K,L\conc(0))$
if the sequences $K$ and $L$ have the same length.

\begin{proposition}\label{proposition:low5:2}
We have
$
\E_j(K,L)\bigl(fH\bigl)=\E_j(K,L)\bigl(f\bigl)\cdot\theta_j^{\smallsum L}(H)
$
for any $f\in\bar\U(n)$ and $H\in\bar\U^0(n)$.
\end{proposition}
\Proof The result follows by induction on $|K|+|L|$ from
Definition~\ref{definition:low5:1} and
Proposition~\ref{proposition:low5:0}\ref{proposition:low5:0:part:2}.
\endproof

For integers $i\le j$ and $k$, we define
$k^{(i,j)}:=\min\{j,\max\{i,k\}\}=$\linebreak$\max\{i,\min\{j,k\}\}$.
For a sequence of integers $K=(k_1,\ldots,k_m)$ and integers $i\le j$,
we put
$
K^{(i,j)}:=(k_1^{(i,j)},\ldots,k_m^{(i,j)})$ and
$
K^{\{j\}}:=(\delta_{k_1\le j},\ldots,\delta_{k_m\le j}).
$\label{minmax}

Below we collect the main properties of the above operations
that we shall use throughout this paper.

\begin{proposition}\label{proposition:low5:3}
Let $K=(k_1,\ldots,k_d)$ be a sequence of integers and $i,j,s,t$ be integers
such that $i\le j$ and $s\le t$.
\begin{enumerate}
\item\label{proposition:low5:3:part:1} If $[i..j]\cap[s..t]\ne\emptyset$
     then $\bigl(K^{(i,j)}\bigl)^{(s,t)}=K^{(\max\{i,s\},\min\{j,t\})}$.
\item\label{proposition:low5:3:part:2} If $i\le s<j$ then
     $\bigl(K^{(i,j)}\bigl)^{\{s\}}=K^{\{s\}}$.
\item\label{proposition:low5:3:part:3}
     $(K-\epsilon_x)^{(i,j)}=K^{(i,j)}-\delta_{i<k_x\le j}\epsilon_x$
     for $1\le x\le d$.
\item\label{proposition:low5:3:part:4} $(K-\epsilon_x)^{\{j\}}=
     K^{\{j\}}+\delta_{k_x-1=j}\epsilon_x$ for $1\le x\le d$.
\end{enumerate}
\end{proposition}

Now we are going to prove a lemma that will explain the role of
sequences $K^{(i,j)}$ and $K^{\{j\}}$.

\begin{lemma}\label{lemma:low5:2}
Suppose $1\le a_1<\cdots<a_q\le n$, where $q\ge 2$, and for each
$i=1,\ldots,q-1$, we are given an element
$X_i\in\bar\U^{-,0}(a_i,a_{i+1})\oplus\bar\I_{M_i}$,
where $M_i\subset[a_i..a_{i+1})$.
Let $K$ be a weakly decreasing sequence of integers with entries from
$[a_1..a_q]$ and $L$ be
a sequence of integers of length $|K|$ or $|K|+1$.
We have
$$
\begin{array}{l}
\E_{a_q}(K,L)(X_1\cdots X_{q-1})\equiv \\[6pt]
\E_{a_2}(K^{(a_1,a_2)},K^{\{a_2\}})(X_1)\cdots\E_{a_{q-1}}(K^{(a_{q-2},a_{q-1})},K^{\{a_{q-1}\}})(X_{q-2})\times\\[6pt]
\times\E_{a_q}(K^{(a_{q-1},a_q)},L)(X_{q-1})\pmod{\bar\I_{M_1\cup\cdots\cup M_{q-1}\cup[\min K..a_q)}}.
\end{array}
$$
\end{lemma}
\Proof
Let $K=(k_1,\ldots,k_m)$ and $L=(l_1,\ldots,l_t)$.
We apply induction on $m+t$. The case $m=t=0$ is obvious.
Therefore assume $m+t>0$.

First we consider the case $m=t$.
Let $\hat K:=(k_1,\ldots,k_{m-1})$ and
$$
\begin{array}{ccl}
\hat X_1  &:=    &\E_{a_2}(\hat K^{(a_1,a_2)},\hat K^{\{a_2\}})(X_1),\\
\vdots    &      &                                                  \\
\hat X_{q-2}&:=    &\E_{a_{q-1}}(\hat K^{(a_{q-2},a_{q-1})},\hat K^{\{a_{q-1}\}})(X_{q-2}),\\[3pt]
\hat X_{q-1}&:=    &\E_{a_q}(\hat K^{(a_{q-1},a_q)}, L)(X_{q-1}).
\end{array}
$$
We put $\hat c:=\min\hat K$ for brevity. Clearly, $k_1=\min K\le\hat c$.
By Definition~\ref{definition:low5:1} and the inductive hypothesis,
we get
$$
\E_{a_q}(K,L)(X_1\cdots X_{q-1})\equiv \E(k_m,a_q)
\hat X_1\cdots\hat X_{q-1}\!
\pmod{\bar\I_{M_1\cup\cdots\cup M_{q-1}\cup[\hat c..a_q)}}.
$$
The required result follows from this formula if $k_m=a_q$.
Therefore, we consider the case $k_m<a_q$.
Let $r=1,\ldots,q-1$ be a number such that
$a_1<\cdots<a_r\le k_m<a_{r+1}<\cdots<a_q$.
Hence
$$
\E(k_m,a_q)=\E(k_m,a_{r+1})\E(a_{r+1},a_{r+2})\cdots\E(a_{q-1},a_q),
$$
where all factors in the right-hand side are nonunitary.
By Proposition~\ref{proposition:low5:-1&low5:2.5}\ref{proposition:low5:2.5}
for every $i=1,\ldots,q-1$, we have
$
\hat X_i\in\bar\U^{-,0}(a_i,a_{i+1})\oplus
\bar\I_{M_i\cup([a_i..a_{i+1})\cap[\hat c..a_q))}
$.
Lemma~\ref{lemma:low5:1}\ref{lemma:low5:1:part:2} now implies
$$
\begin{array}{l}
\E(k_m,a_q)\hat X_1\cdots\hat X_{q-1}\equiv
\hat X_1\cdots\hat X_{r-1}\E(k_m,a_{r+1})\times\\
\times\theta_{a_{r+1}}(\hat X_r)\E(a_{r+1},a_{r+2})\cdots
      \theta_{a_{q-1}}(\hat X_{q-2})\E(a_{q-1},a_q)\times\\
\times\hat X_{q-1}\pmod{\bar\I_{M_1\cup\cdots\cup M_{q-1}\cup[k_1..a_q)}}
\end{array}
$$
if $r=1$ or $a_r<k_m$ and
$$
\begin{array}{l}
\E(k_m,a_q)\hat X_1\cdots\hat X_{q-1}\equiv
\hat X_1\cdots\hat X_{r-2}\theta_{a_r}(\hat X_{r-1})\E(k_m,a_{r+1})\times\\
\times\theta_{a_{r+1}}(\hat X_r)\E(a_{r+1},a_{r+2})\cdots
      \theta_{a_{q-1}}(\hat X_{q-2})\E(a_{q-1},a_q)\times\\
\times\hat X_{q-1}\pmod{\bar\I_{M_1\cup\cdots\cup M_{q-1}\cup[k_1..a_q)}}
\end{array}
$$
if $r>1$ and $a_r=k_m$.
The result now follows from Definition~\ref{definition:low5:1}.

Now let $t=m+1$. By the inductive hypothesis and
Proposition~\ref{proposition:low5:0}\ref{proposition:low5:0:part:5}, we get
\begin{equation}\label{equation:low:5:7.5}
\begin{array}{l}
\E_{a_q}(K,L)(\X_1\cdots\X_{q-1})\equiv \theta_{a_q}^{l_t}
\Bigl(
\E_{a_2}(K^{(a_1,a_2)},K^{\{a_2\}})(X_1)\cdots\\[6pt]
\E_{a_{q-1}}(K^{(a_{q-2},a_{q-1})},K^{\{a_{q-1}\}})(X_{q-2})
\E_{a_q}(K^{(a_{q-1},a_q)},\hat L)(X_{q-1})
\Bigr)
\\[6pt]
\pmod{\bar\I_{M_1\cup\cdots\cup M_{q-1}\cup[\min K..a_q)}},
\end{array}
\end{equation}
where $\hat L=(l_1,\ldots,l_m)$.
It follows from Proposition~\ref{proposition:low5:-1&low5:2.5}\ref{proposition:low5:2.5}
that for every $i=1,\ldots,q-1$, there exists
$\hat Y_i\in\bar\U^{-,0}(a_i,a_{i+1})$ such that
$$
\E_{a_{i+1}}(K^{(a_i,a_{i+1})},N_i)(X_i)
\equiv \hat Y_i\pmod{\bar\I_{M_i\cup([a_i..a_{i+1})\cap[\min K..a_q))}},
$$
where $N_i=K^{\{a_{i+1}\}}$ if $i<q-1$ and $N_{q-1}=\hat L$.
It follows from Lemma~\ref{lemma:low5:1}\ref{lemma:low5:1:part:1},
Proposition~\ref{proposition:low5:0}\ref{proposition:low5:0:part:5},\ref{proposition:low5:0:part:3}
and formula~(\ref{equation:low:5:7.5}) that
$$
\begin{array}{l}
\E_{a_q}(K,L)(X_1\cdots X_{q-1})\equiv
\hat Y_1\cdots\hat Y_{q-2}\theta_{a_q}^{l_t}(\hat Y_{q-1})\\[6pt]
\equiv \E_{a_2}(K^{(a_1,a_2)},K^{\{a_2\}})(X_1)\cdots\E_{a_{q-1}}(K^{(a_{q-2},a_{q-1})},K^{\{a_{q-1}\}})(X_{q-2})\times\\[6pt]
\times\theta_{a_q}^{l_t}(\E_{a_q}(K^{(a_{q-1},a_q)},\hat L)(X_{q-1}))\pmod{\bar\I_{M_1\cup\cdots\cup M_{q-1}\cup[\min K..a_q)}}.
\end{array}
$$
It remains to apply Definition~\ref{definition:low5:1} to the last factor.
\endproof

\section{Multiplication formulae}\label{Multiplication formulae}

In this section, we introduce certain elements $\T_{i,j}^{(d)}(M,R)$
of $\bar\U(n)$. We are interested in the behavior of
$\E_j(K,L)\bigl(\T_{i,j}^{(d)}(M,R)\bigr)$ under
the left multiplication by $\E_{l-1}$
modulo a suitable ideal. This is done in Lemmas~\ref{lemma:low5:3}
and~\ref{lemma:low5:4}. In what follows, $\Z(\zeta)$ denotes
the field of fractions of the ring of polynomials $\Z[\zeta]$.

\begin{lemma}\label{lemma:low5:2.5}
Let $1\le i<j\le n$, $d\ge1$, $R\in\Z(\zeta)$
and $K=(k_1,\ldots,k_m)$, $L=(l_1,\ldots,l_q)$ be sequences of integers
such that $m\le d$, $i\le k_1\le\cdots\le k_m\le j$ and $q=m$ or $q=m+1$.
Modulo $\bar\I_{[\min K'..j)}$, we have
$$
\begin{array}{l}
\displaystyle \E_j(K,L)\Bigl(\S_{i,j}^dR(\C(i,j))\Bigr)\equiv \S_{i,K'}\;d^{\underline r}\times\\[6pt]
\displaystyle \times \prod_{s=1}^r\left(
\left(\C(i,j)-d+s-1+\sum_{s<h\le q}l_h\right)\prod_{k_s<t<j}(\C(i,t)-d+s-1)
\right)\times\\[22pt]
\times R(\C(i,j)+\smallsum L),
\end{array}
$$
where $K'=(k_1,\ldots,k_m,j^{d-m})$ and
$r$ is the number of entries of $K$ less than $j$.
\end{lemma}
\Proof We apply induction on $m+q$. If $q=m+1$ then
we obtain the required formula from the inductive hypothesis
by applying Definition~\ref{definition:low5:1} and
Proposition~\ref{proposition:low5:0}\ref{proposition:low5:0:part:4},\ref{proposition:low5:0:part:5}.

Now let $m=q$.
Since the case $K=\emptyset$ is obvious, we assume $K\ne\emptyset$.
By definition~\ref{definition:low5:1} we have
$$
\E_j(K,L)\Bigl( \S_{i,j}^dR(\C(i,j)) \Bigr)=
\E(k_m,j)\E_j(\hat K,L)\Bigl(\S_{i,j}^dR(\C(i,j))\Bigr),
$$
where $\hat K=(k_1,\ldots,k_{m-1})$.
The inductive hypothesis implies
{\arraycolsep=-1pt
\begin{equation}\label{equation:low5:8}
\!\!
\begin{array}{l}
\displaystyle \E_j(\hat K,L)\Bigl(\S_{i,j}^dR(\C(i,j))\Bigl)\equiv \S_{i,K''}d^{\underline{\hat r}}\times\\[6pt]
\displaystyle
\times\prod_{s=1}^{\hat r}\left(
\left(\C(i,j)-d+s-1+\sum_{s<h\le m}l_h\right)\prod_{k_s<t<j}(\C(i,t)-d+s-1)
\right)\\[22pt]
\times R(\C(i,j)+\smallsum L)\pmod{\bar\I_{[\min K'..j)}},
\end{array}
\end{equation}}
\!\!where $K''=(k_1,\ldots,k_{m-1},j^{d-m+1})$ and
$\hat r$ is the number of entries of $\hat K$ less than $j$.
Clearly, $\min K''\ge\min K'$.
If $k_m=j$ then $\E(k_m,j)=1$, $\hat r=r$ and $K''=K'$.
Hence the required result.

Now let $k_m<j$. We have $r=m$ and $\hat r=m-1$.
By Corollary~\ref{corollary:low3:1} we obtain
$$
\E(k_m,j)\S_{i,K''}\equiv \S_{i,K'}(d-m+1)\prod_{k_m<t\le j}(\C(i,t)-(d-m+1))
\pmod{\bar\I_{[k_m..j)}}.
$$
Using this formula and multiplying~(\ref{equation:low5:8}) by $\E(k_m,j)$
on the left, we obtain the required formula. \endproof

Now we introduce the central object of our study.

\begin{definition}\label{definition:low5:2}
Let $1\le i<j\le n$, $d\ge1$ and $R\in\Z(\zeta)$
We define
\begin{itemize}
\item $\T_{i,j}^{(d)}(\emptyset,R)=\S_{i,j}^dR(\C(i,j))$;
\item $\T_{i,j}^{(d)}(M,R)=\bigl(\T_{i,j}^{(d)}(M',R)-\S_{i,m}^d
       \T_{m,j}^{(d)}(M',R)\bigl)\C(i,m)^{-1}$,
       where $m=\min M$ and $M'=M\setminus\{m\}$,
       if $M$ is a nonempty subset of $(i..j)$.
\end{itemize}
\end{definition}
The elements $\T_{i,j}^{(d)}(M,R)$ belong to $\bar\U^{-,0}(i,j)$.

\begin{lemma}\label{lemma:low5:3}
Let $1\le i<j\le n$, $d\ge1$, $M\subset(i..j)$, $R\in\Z(\zeta)$,
$K=(k_1,\ldots,k_d)$ and $L$ be sequences of integers such that
$i\le k_1\le\cdots\le k_d\le j$ and $|L|=d$ or $|L|=d+1$.
Take $l\in(1..n]\setminus K$.
Then modulo the ideal $\bar\I_{[i..j)\cup\{l-1\}}$,
and even the ideal $\bar\I_{\{l-1\}}$ if $K=(j^d)$, we have
$$
\E_{l-1}\E_j(K,L)\bigl(\T_{i,j}^{(d)}(M,R)\bigl)\equiv 0
$$
if $l\notin M$ or  $k_d<l$ and
$$
\begin{array}{l}
\E_{l-1}\E_j(K,L)\bigl(\T_{i,j}^{(d)}(M,R)\bigl)\equiv \\[6pt]
-\left[\E_l(K^{(i,l)}-\epsilon_{a+1},K^{\{l\}})\left(\T_{i,l}^{(d)}\left(M\cap(i..l),\frac1{\zeta-d}\right)\right)\right]\times \\[6pt]
\times\E_j(K^{(l,j)},L)\bigl(\T_{l,j}^{(d)}(M\cap(l..j),R)\bigl)
\end{array}
$$
otherwise, where $a$ is the number of entries of $K$ less than $l$.
\end{lemma}
\Proof
Let $P:=\emptyset$ if $K=(j^d)$ and $P:=[i..j)$ otherwise.
We apply induction on $|M|$.
Let $M=\emptyset$.
By Lemma~\ref{lemma:low5:2.5}, we have
$\E_j(K,L)\bigl(\T_{i,j}^{(d)}(\emptyset,R)\bigl)\equiv
\S_{i,K}H\pmod{\bar\I_P}$ for some $H\in\bar\U^0(n)$.
Now the result follows from~(\ref{equation:low2:0.5}) and $l\notin K$.

Now let $M\ne\emptyset$ and the result holds for smaller sets.
We put $m:=\min M$ and $M':=M\setminus\{m\}$.
By Definition~\ref{definition:low5:2}, Proposition~\ref{proposition:low5:2}
and Lemma~\ref{lemma:low5:2} we have
\begin{equation}\label{equation:low5:2.5}
\begin{array}{l}
\E_j(K,L)\bigl(\T_{i,j}^{(d)}(M,R)\bigl)\equiv
\Bigl(
\E_j(K,L)\bigl(\T_{i,j}^{(d)}(M',R)\bigl)-\\[6pt]
\E_m(K^{(i,m)},K^{\{m\}})\bigl(\S_{i,m}^d\bigl)\cdot
\E_j(K^{(m,j)},L)\bigl(\T_{m,j}^{(d)}(M',R)\bigl)
\Bigl)\C(i,m)^{-1}\\[6pt]
\pmod{\bar\I_P}.
\end{array}
\end{equation}

We are going to multiply this equivalence by $\E_{l-1}$ on the left and apply
the inductive hypothesis.

{\it Case $\it l\notin M$ or $k_d<l$}. The inductive hypothesis implies
\begin{equation}\label{equation:low5:3.5}
\E_{l-1}\E_j(K,L)\bigl(\T_{i,j}^{(d)}(M',R)\bigl)\equiv 0\pmod{\bar\I_{P\cup\{l-1\}}}.
\end{equation}
If $l-1\notin[i..j)$ then
by Proposition~\ref{proposition:low5:-1&low5:2.5}\ref{proposition:low5:2.5}
and Lemma~\ref{lemma:low5:1}\ref{lemma:low5:1:part:1} we have
\begin{equation}\label{equation:low5:4}
\begin{array}{l}
\E_{l-1}\E_m(K^{(i,m)},K^{\{m\}})\bigl(\S_{i,m}^d\bigl)\times\\[6pt]
\times
\E_j(K^{(m,j)},L)\bigl(\T_{m,j}^{(d)}(M',R)\bigl)\equiv 0\pmod{\bar\I_{P\cup\{l-1\}}}.
\end{array}
\end{equation}
If $l-1\in[i..m)$ then $l\notin K^{(i,m)}$.
Therefore, the inductive hypothesis implies
$\E_{l-1}\E_m(K^{(i,m)},K^{\{m\}})\bigl(\S_{i,m}^d\bigl)\equiv 0$ modulo
$\bar\I_{(P\cap [i..m))\cup\{l-1\}}$. If we multiply the last ideal
by $\E_j(K^{(m,j)},L)\bigl(\T_{m,j}^{(d)}(M',R)\bigl)$ on the right,
then we obtain a subset of $\bar\I_{P\cup\{l-1\}}$.
Hence we get~(\ref{equation:low5:4}) again.

Finally assume $l-1\in[m..j)$.
By proposition~\ref{proposition:low5:-1&low5:2.5}\ref{proposition:low5:2.5}
and Lemma~\ref{lemma:low5:1}\ref{lemma:low5:1:part:2}
we get
\begin{equation}\label{equation:low5:5}
\!\!\!
\begin{array}{l}
\E_{l-1}\E_m(K^{(i,m)},K^{\{m\}})\bigl(\S_{i,m}^d\bigl)\cdot\E_j(K^{(m,j)},L)\bigl(\T_{m,j}^{(d)}(M',R)\bigl)\equiv \\[6pt]
       \theta_m^{\delta_{l-1=m}}\Bigl(\E_m(K^{(i,m)},K^{\{m\}})\bigl(\S_{i,m}^d\bigl)\Bigl)\cdot\E_{l-1}\E_j(K^{(m,j)},L)\bigl(\T_{m,j}^{(d)}(M',R)\bigl)\\[6pt]
\pmod{\bar\I_{P\cup\{l-1\}}}.
\end{array}
\end{equation}
Clearly $l\notin K^{(m,j)}$ and all entries of $K^{(m,j)}$ are less than $l$
if $k_d<l$. Therefore, by the inductive hypothesis
$\E_{l-1}\E_j(K^{(m,j)},L)\bigl(\T_{m,j}^{(d)}(M',R)\bigl)\equiv 0$ modulo $\bar\I_{P\cup\{l-1\}}$,
whence we have~(\ref{equation:low5:4}) one more time.

{\it Case $l=m$ and $k_d\ge l$}. By the inductive hypothesis
formula~(\ref{equation:low5:3.5}) holds.
From Definition~\ref{definition:low5:1}, it follows that
$$
\E_{l-1}\E_m(K^{(i,m)},K^{\{m\}})\bigl(\S_{i,m}^d\bigl)=
\E_m(K^{(i,m)}-\epsilon_{a+1},K^{\{m\}})\bigl(\S_{i,m}^d\bigl).
$$
It remains to deal with $\C(i,m)^{-1}$. Firstly, we have
$$
\begin{array}{l}
\E_j(K^{(m,j)},L)\bigl(\T_{m,j}^{(d)}(M',R)\bigr)\cdot\C(i,m)^{-1}=\\[6pt]
(\C(i,m)-d+a)^{-1}\cdot\E_j(K^{(m,j)},L)\bigl(\T_{m,j}^{(d)}(M',R)\bigr).
\end{array}
$$
Secondly, applying Proposition~\ref{proposition:low5:2} and
Proposition~\ref{proposition:low5:0}\ref{proposition:low5:0:part:4}, we get
$$
\begin{array}{l}
\E_m(K^{(i,m)}-\epsilon_{a+1},K^{\{m\}})\bigl(\S_{i,m}^d\bigl)\cdot
(\C(i,m)-d+a)^{-1}\\[6pt]
=\E_m(K^{(i,m)}-\epsilon_{a+1},K^{\{m\}})\bigl(\S_{i,m}^d(\C(i,m)-d)^{-1}\bigl).
\end{array}
$$

{\it Case $l\in M'$ and $k_d\ge l$}.
Multiplying~(\ref{equation:low5:2.5}) by $\E_{l-1}$ on the left and
applying Proposition~\ref{proposition:low5:3} and Lemma~\ref{lemma:low5:1},
we obtain
\begin{equation}\label{equation:low5:6}
\begin{array}{l}
\E_{l-1}\E_j(K,L)\bigl(\T_{i,j}^{(d)}(M,R)\bigl)\equiv \\[6pt]
-\Biggl(
\E_l(K^{(i,l)}-\epsilon_{a+1},K^{\{l\}})\left(\T_{i,l}^{(d)}\left(M'\cap(i..l),\frac1{\zeta-d}\right)\right)\\[16pt]
-\biggl[\theta_m^{\delta_{l-1=m}}\Bigl(\E_m(K^{(i,m)},K^{\{m\}})\bigl(\S_{i,m}^d\bigl)\Bigl)\times \\[6pt]
\times\E_l(K^{(m,l)}-\epsilon_{a+1},K^{\{l\}})\left(\T_{m,l}^{(d)}\left(M'\cap(m..l),\frac1{\zeta-d}\right)\right)\biggr]
\Biggr)\times \\[12pt]
\times\E_j(K^{(l,j)},L)\bigl(\T_{l,j}^{(d)}(M\cap(l..j),R))\C(i,m)^{-1}
\pmod{\bar\I_{P\cup\{l-1\}}}.
\end{array}
\end{equation}
Applying Lemma~\ref{lemma:low5:2} and Proposition~\ref{proposition:low5:3},
we get
\begin{equation}\label{equation:low5:7}
\begin{array}{l}
\E_l(K^{(i,l)}-\epsilon_{a+1},K^{\{l\}})
\left(\S_{i,m}^d\T^{(d)}_{m,l}\left(M'\cap(m..l),\frac1{\zeta-d}\right)\right)\equiv \\[6pt]
\E_m(K^{(i,m)},K^{\{m\}}+\delta_{l-1=m}\epsilon_{a+1})\bigl(\S_{i,m}^d\bigl)\times\\[6pt]
\times\E_l(K^{(m,l)}-\epsilon_{a+1},K^{\{l\}})
\left(\T^{(d)}_{m,l}\left(M'\cap(m..l),\frac1{\zeta-d}\right)\right)\\[6pt]
\!\!\!\!\!\pmod{\bar\I_{(P\cap[i..l))\cup\{l-1\}}}.
\end{array}
\end{equation}
Moreover, Proposition~\ref{proposition:low5:-1&low5:2.5}\ref{proposition:low5:2.5}
shows that both sides of~(\ref{equation:low5:7}) belong to\linebreak
$\bar\U^{-,0}(i,l)\oplus\bar\I_{(P\cap [i..l))\cup\{l-1\}}$.
Finally, we observe
\begin{equation}\label{equation:low5:8.5}
\E_m(K^{(i,m)},K^{\{m\}}+\delta_{l-1=m}\epsilon_{a+1})=
\theta_m^{\delta_{l-1=m}}\circ\E_m(K^{(i,m)},K^{\{m\}}).
\end{equation}
Indeed, in the only nontrivial case $l-1=m$,
this formula follows from Definition~\ref{definition:low5:1}
and the equalities $K^{(i,m)}=(k_1,\ldots,k_a,m^{d-a})$ and
$K^{\{m\}}=(1^a,0^{d-a})$.
Therefore, the right-hand side of~(\ref{equation:low5:7})
is exactly the product in the
square brackets of~(\ref{equation:low5:6}).
Applying Lemma~\ref{lemma:low5:1}\ref{lemma:low5:1:part:1}, we get
$$
\begin{array}{l}
\E_{l-1}\E_j(K,L)\bigl(\T_{i,j}^{(d)}(M,R)\bigl)\equiv \\[6pt]
-\E_l(K^{(i,l)}-\epsilon_{a+1},K^{\{l\}})
\Bigl(
\T_{i,l}^{(d)}\left(M'\cap(i..l),\frac1{\zeta-d}\right)-\\[6pt]
\S_{i,m}^d\T^{(d)}_{m,l}\left(M'\cap(m..l),\frac1{\zeta-d}\right)
\Bigr)\C(i,m)^{-1}\times\\[6pt]
\times\E_j(K^{(l,j)},L)\bigl(\T_{l,j}^{(d)}(M\cap(l..j),R))
\pmod{\bar\I_{P\cup\{l-1\}}}.
\end{array}
$$
Applying Proposition~\ref{proposition:low5:2},
we carry $\C(i,m)^{-1}$ (unchanged) under
the sign of $\E_l(K^{(i,l)}-\epsilon_{a+1},K^{\{l\}})$
and use Definition~\ref{definition:low5:2}.
\endproof

\begin{lemma}\label{lemma:low5:4}
Let $1\le i<j\le n$, $d\ge1$, $M\subset(i..j)$, $R\in\Z(\zeta)$
and $K=(k_1,\ldots,k_d)$, $L=(l_1,\ldots,l_q)$ be
sequences of integers such that $i\le k_1\le\cdots\le k_d\le j$
and $q=d$ or $q=d+1$. Take $l\in(i..j]\cap K$.
Then modulo $\bar\I_{[i..j)}$, or even modulo the zero ideal if
$l=j$, we have
$$
\E_{l-1}\E_j(K,L)\bigl(\T_{i,j}^{(d)}(M,R)\bigl)\equiv
c\,\E_j(K-\epsilon_{a+1},L')\bigl(\T_{i,j}^{(d)}(M,R)\bigl)
$$
if $l\notin M$ and
$$
\begin{array}{l}
\E_{l-1}\E_j(K,L)\bigl(\T_{i,j}^{(d)}(M,R)\bigl)\equiv
c\,\E_j(K-\epsilon_{a+1},L)\bigl(\T_{i,j}^{(d)}(M,R)\bigl)+\\[6pt]
\left[(c-1)\,\E_l(K^{(i,l)}-\epsilon_{a+1},K^{\{l\}})\left(\T_{i,l}^{(d)}\left(M\cap(i..l),\frac1{\zeta-d-1}\right)\right)\right]\times\\[6pt]
\times\E_j(K^{(l,j)},L)\bigl(\T_{l,j}^{(d)}(M\cap(l..j),R)\bigl)
\end{array}
$$
if $l\in M$, where $a$ is the number of entries of $K$ less than $l$;
$L'=L$ if $l<j$ and $L'=(l_1,\ldots,l_a,\sum_{a<h\le q}l_h,0^{q-a-1})$ if $l=j$;
$c$ is the number of entries of $K$ equal to $l$ if $l<j$ and $c=1$ if $l=j$.
\end{lemma}
\Proof If $l=j$ then Definition~\ref{definition:low5:1} implies
$$
\E_j(K,L)=\theta_j^{{\textstyle\Sigma}_{a<h\le q}l_h}\circ\E_j\bigl((k_1,\ldots,k_a),(l_1,\ldots,l_a)\bigl),
$$
whence the required result easily follows. Therefore we assume $l<j$
for the rest of the proof. Let $b$ and $r$ denote the number
of entries of $K$ not greater than $l$ and less than $r$ respectively.
We have $c=b-a$ and $L'=L$. Since $l<j$, the number of entries of
$K-\epsilon_{a+1}$ less than $j$ is also $r$.
We apply induction on $|M|$.
In the case $M=\emptyset$, the required result immediately follows from
Lemma~\ref{lemma:low5:2.5} and Corollary~\ref{corollary:low3:1}.

Now let $M\ne\emptyset$ and suppose that the lemma is true for
sets of lesser cardinality.
We put $m:=\min M$ and $M':=M\setminus\{m\}$.
Clearly,~(\ref{equation:low5:2.5}) holds in the present case,
where $P=[i..j)$.

{\it Case $\it l\notin M$}. By the inductive hypothesis, we have
\begin{equation}\label{equation:low5:10.5}
\begin{array}{l}
\E_{l-1}\E_j(K,L)\bigl(\T_{i,j}^{(d)}(M',R)\bigl)\\[6pt]
\equiv (b-a)\E_j(K-\epsilon_{a+1},L)\bigl(\T_{i,j}^{(d)}(M',R)\bigl)
\pmod{\bar\I_{[i..j)}}.
\end{array}
\end{equation}

We first consider the case $i<l<m$.
The inductive hypothesis and Proposition~\ref{proposition:low5:3} yield
$$
\begin{array}{l}
\E_{l-1}\E_m(K^{(i,m)},K^{\{m\}})\bigl(\S_{i,m}^d\bigl)\equiv \\[6pt]
(b-a)\E_m((K-\epsilon_{a+1})^{(i,m)},(K-\epsilon_{a+1})^{\{m\}})\bigl(\S_{i,m}^d\bigl)
\pmod{\bar\I_{[i..m)}}.
\end{array}
$$
Applying Proposition~\ref{proposition:low5:-1&low5:2.5}\ref{proposition:low5:2.5},
Lemma~\ref{lemma:low5:1}\ref{lemma:low5:1:part:1} and
Proposition~\ref{proposition:low5:3}\ref{proposition:low5:3:part:3},
we get
\begin{equation}\label{equation:low5:11}
\begin{array}{l}
\E_{l-1}\E_m(K^{(i,m)},K^{\{m\}})\bigl(\S_{i,m}^d\bigl)\cdot\E_j(K^{(m,j)},L)\bigl(\T_{m,j}^{(d)}(M',R)\bigl)\\[6pt]
\equiv (b-a)\E_m((K-\epsilon_{a+1})^{(i,m)},(K-\epsilon_{a+1})^{\{m\}})\bigl(\S_{i,m}^d\bigl)\times\\[6pt]
\times\E_j((K-\epsilon_{a+1})^{(m,j)},L)\bigl(\T_{m,j}^{(d)}(M',R)\bigl)
\pmod{\bar\I_{[i..j)}}.
\end{array}
\end{equation}

Now we consider the case $m<l<j$.
In that case,~(\ref{equation:low5:5}) also holds, where $P=[i..j)$
(i.e., the equivalence holds modulo $\bar\I_{[i..j)}$).
To rearrange the second factor of
the right-hand side of~(\ref{equation:low5:5}), we apply
the inductive hypothesis and Proposition~\ref{proposition:low5:3}. Thus we get
$$
\begin{array}{l}
\E_{l-1}\E_j(K^{(m,j)},L)\bigl(\T_{m,j}^{(d)}(M',R)\bigl)\equiv \\[6pt]
(b-a)\E_j((K-\epsilon_{a+1})^{(m,j)},L)\bigl(\T_{m,j}^{(d)}(M',R)\bigl)
\pmod{\bar\I_{[m..j)}}.
\end{array}
$$
Since in the present case~(\ref{equation:low5:8.5}) also holds,
we can apply this formula to the first factor of
the right-hand side of~(\ref{equation:low5:5}).
Finally, Proposition~\ref{proposition:low5:3} shows
that~(\ref{equation:low5:11}) is again true.

In both cases,~(\ref{equation:low5:11}) and Lemma~\ref{lemma:low5:2} yield
\begin{equation}\label{equation:low5:11.5}
\begin{array}{l}
\E_{l-1}\E_m(K^{(i,m)},K^{\{m\}})\bigl(\S_{i,m}^d\bigl)\cdot\E_j(K^{(m,j)},L)\bigl(\T_{m,j}^{(d)}(M',R)\bigl)\\[6pt]
\equiv (b-a)\E_j(K-\epsilon_{a+1},L)\bigl(\S_{i,m}^d\T_{m,j}^{(d)}(M',R)\bigl)
\pmod{\bar\I_{[i..j)}}.
\end{array}
\end{equation}

Therefore, multiplying~(\ref{equation:low5:2.5}) on the left by $\E_{l-1}$,
using~(\ref{equation:low5:10.5}) and~(\ref{equation:low5:11.5}),
Proposition~\ref{proposition:low5:2} to carry $\C(i,m)^{-1}$ under
the sign of $\E_j(K-\epsilon_{a+1},L)$,
and Definition~\ref{definition:low5:2}, we obtain the required result.

{\it Case $l=m$}. The inductive hypothesis implies
that~(\ref{equation:low5:10.5}) holds.
We introduce the auxiliary notation $K':=(1^a,b-a,0^{d-a-1})$ and
$$
\begin{array}{l}
X:=(b-a)\E_m((K-\epsilon_{a+1})^{(i,m)},(K-\epsilon_{a+1})^{\{m\}})
\bigl(\S_{i,m}^d\bigl)\\[6pt]
-\E_m(K^{(i,m)}-\epsilon_{a+1},K')\bigl(\S_{i,m}^d\bigl).
\end{array}
$$
It follows from
Proposition~\ref{proposition:low5:-1&low5:2.5}\ref{proposition:low5:2.5} that
$X\in\bar\U^{-,0}(i,m)\oplus\bar\I_{[i..m)}$.
With regard to
Proposition~\ref{proposition:low5:3}\ref{proposition:low5:3:part:2},
the inductive hypothesis implies
\begin{equation}\label{equation:low5:11.75}
\E_{l-1}\E_m(K^{(i,m)},K^{\{m\}})\bigl(\S_{i,m}^d\bigl)
=\E_m(K^{(i,m)}-\epsilon_{a+1},K')\bigl(\S_{i,m}^d\bigl).
\end{equation}

Therefore, multiplying~(\ref{equation:low5:2.5}) by $\E_{l-1}$
on the left, using~(\ref{equation:low5:10.5}) and~(\ref{equation:low5:11.75}),
and applying Proposition~\ref{proposition:low5:3}\ref{proposition:low5:3:part:3},
we get
$$
\begin{array}{l}
\E_{l-1}\E_j(K,L)\bigl(\T_{i,j}^{(d)}(M,R)\bigl)\equiv
\Bigl(
(b-a)\E_j(K-\epsilon_{a+1},L)\bigl(\T_{i,j}^{(d)}(M',R)\bigl)-\\[6pt]
\E_m(K^{(i,m)}-\epsilon_{a+1},K')\bigl(\S_{i,m}^d\bigl)\cdot
\E_j(K^{(m,j)},L)\bigl(\T_{m,j}^{(d)}(M',R)\bigl)
\Bigl)\C(i,m)^{-1}=\\[6pt]
(b-a)\biggl(
\E_j(K-\epsilon_{a+1},L)\bigl(\T_{i,j}^{(d)}(M',R)\bigl)-\\[6pt]
\Bigl[\E_m((K-\epsilon_{a+1})^{(i,m)},(K-\epsilon_{a+1})^{\{m\}})\bigl(\S_{i,m}^d\bigl)\times\\[6pt]
\times\E_j((K-\epsilon_{a+1})^{(m,j)},L)\bigl(\T_{m,j}^{(d)}(M',R)\bigl)\Bigr]
\biggl)\C(i,m)^{-1}\\[6pt]
+X\cdot\E_j(K^{(m,j)},L)\bigl(\T_{m,j}^{(d)}(M',R)\bigl)\C(i,m)^{-1}
\pmod{\bar\I_{[i..j)}}
\end{array}
$$
By Lemma~\ref{lemma:low5:2} the product in the square brackets of
the above formula equals
$
\E_j(K-\epsilon_{a+1},L)\bigl(\S_{i,m}^d\T_{m,j}^{(d)}(M',R)\bigl)
$
modulo $\bar\I_{[i..j)}$.
Thus we obtain
\begin{equation}\label{equation:low5:12}
\begin{array}{l}
\E_{l-1}\E_j(K,L)\bigl(\T_{i,j}^{(d)}(M,R)\bigl)
\equiv (b-a)\E_j(K-\epsilon_{a+1},L)\bigl(\T_{i,j}^{(d)}(M,R)\bigl)\\[6pt]
+X\cdot\E_j(K^{(m,j)},L)\bigl(\T_{m,j}^{(d)}(M',R)\bigl)\C(i,m)^{-1}
\pmod{\bar\I_{[i..j)}}.
\end{array}
\end{equation}

Let us calculate $X$ modulo $\bar\I_{[i..m)}$.
By Lemma~\ref{lemma:low5:2.5}, we have
{\arraycolsep=1pt
\begin{equation}\label{equation:low5:12.75}
\begin{array}{l}
\displaystyle X\equiv \S_{i,K^{(i,m)}-\epsilon_{a+1}}d^{\underline{a+1}}
\left(\prod_{s=1}^a\prod_{k_s<t<m}(\C(i,t)-d+s-1)\right)\times\\[18pt]
\displaystyle
\times
(\C(i,m){-}d{+}b{-}1)^a
\bigl((b{-}a)(\C(i,m){-}d{+}b{-}1){-}(\C(i,m){-}d{+}a)\bigl)\equiv
\\[12pt]
(b{-}a{-}1)\E_m(K^{(i,m)}{-}\epsilon_{a+1},K^{\{m\}})\bigl(\S_{i,m}^d\bigl)
\frac{\C(i,m)-d+b}{\C(i,m)-d+b-1}
\!\!\pmod{\bar\I_{[i..m)}}.
\end{array}
\end{equation}
}

To rearrange the second summand of the right-hand side
of~(\ref{equation:low5:12}), we apply
Proposition~\ref{proposition:low5:-1&low5:2.5}\ref{proposition:low5:2.5},
equivalence~(\ref{equation:low5:12.75}),
Lemma~\ref{lemma:low5:1}\ref{lemma:low5:1:part:1}
and Propositions~\ref{proposition:low5:2} and~\ref{proposition:low5:0}\ref{proposition:low5:0:part:4}.
Thus we obtain
$$
\begin{array}{l}
X\cdot\E_j(K^{(m,j)},L)\bigl(\T_{m,j}^{(d)}(M',R)\bigl)\C(i,m)^{-1}\\[6pt]
=X(\C(i,m)-d+b)^{-1}\E_j(K^{(m,j)},L)\bigl(\T_{m,j}^{(d)}(M',R)\bigl)\\[6pt]
\equiv (b-a-1)\E_m(K^{(i,m)}-\epsilon_{a+1},K^{\{m\}})\bigl(\S_{i,m}^d\bigl)
(\C(i,m)-d+b-1)^{-1}\times \\[6pt]
\times\E_j(K^{(m,j)},L)\bigl(\T_{m,j}^{(d)}(M',R)\bigl)\\[6pt]
=(b-a-1)\E_m(K^{(i,m)}-\epsilon_{a+1},K^{\{m\}})
\left(\T^{(d)}_{i,m}\left(\emptyset,\frac1{\zeta-d-1}\right)\right)\times\\[6pt]
\times\E_j(K^{(m,j)},L)\bigl(\T_{m,j}^{(d)}(M',R)\bigr)
\pmod{\bar\I_{[i..j)}}.
\end{array}
$$
This formula together with~(\ref{equation:low5:12}) yield
the required result.

{\it Case $l\in M'$.} By the inductive hypothesis, we have
$$
\begin{array}{l}
\E_{l-1}\E_j(K,L)\bigl(\T_{i,j}^{(d)}(M',R)\bigl)\equiv
(b-a)\E_j(K-\epsilon_{a+1},L)\bigl(\T_{i,j}^{(d)}(M',R)\bigl)+\\[6pt]
(b-a-1)\E_l(K^{(i,l)}-\epsilon_{a+1},K^{\{l\}})\left(\T_{i,l}^{(d)}\left(M'\cap(i..l),\frac1{\zeta-d-1}\right)\right)\times\\[6pt]
\times\E_j(K^{(l,j)},L)\bigl(\T_{l,j}^{(d)}(M\cap(l..j),R)\bigl)
\pmod{\bar\I_{[i..j)}}.
\end{array}
$$
In the present case,~(\ref{equation:low5:5}) also holds, where $P=[i..j)$.
To rearrange the second factor of
the right-hand side of~(\ref{equation:low5:5}), we apply
the inductive hypothesis and Proposition~\ref{proposition:low5:3}.
Thus we get
$$
{
\arraycolsep=0pt
\begin{array}{l}
 \E_{l-1}\E_j(K^{(m,j)},L)\bigl(\T_{m,j}^{(d)}(M',R)\bigl)\equiv
 (b{-}a)\E_j(K^{(m,j)}{-}\epsilon_{a+1},L)\bigl(\T_{m,j}^{(d)}(M',R)\bigl)\\[6pt]
 +(b-a-1)\E_l(K^{(m,l)}-\epsilon_{a+1},K^{\{l\}})\left(\T_{m,l}^{(d)}\left(M'\cap(m..l),\frac1{\zeta-d-1}\right)\right)\times\\[6pt]
 \times\E_j(K^{(l,j)},L)\bigl(\T_{l,j}^{(d)}(M\cap(l..j),R)\bigl)
 \pmod{\bar\I_{[m..j)}}.
\end{array}
}
$$
Multiplying~(\ref{equation:low5:2.5}) by $\E_{l-1}$ on the left
and taking into account~(\ref{equation:low5:5}) together
with the above equivalences, we obtain
$$
\begin{array}{l}
\E_{l-1}\E_j(K,L)\bigl(\T_{i,j}^{(d)}(M,R)\bigl)\equiv
(b-a)\biggl(\E_j(K-\epsilon_{a+1},L)\bigl(\T_{i,j}^{(d)}(M',R)\bigl)-\\[6pt]
\Bigl[\theta_m^{\delta_{l-1=m}}\bigl(\E_m(K^{(i,m)},K^{\{m\}})(\S_{i,m}^d)\bigl)
\E_j(K^{(m,j)}{-}\epsilon_{a+1},L)\bigl(\T_{m,j}^{(d)}(M',R)\bigl)\Bigr]\biggr){\times}
\\[6pt]\times
\C(i,m)^{-1}
{+}(b{-}a{-}1)\Biggl(\E_l(K^{(i,l)}{-}\epsilon_{a+1},K^{\{l\}})\left(\T_{i,l}^{(d)}\left(M'\cap(i..l),\frac1{\zeta-d-1}\right)\right)\\[6pt]
-\biggl[\theta_m^{\delta_{l-1=m}}\bigl(\E_m(K^{(i,m)},K^{\{m\}})(\S_{i,m}^d)\bigl)\times\\[6pt]
\times\E_l(K^{(m,l)}-\epsilon_{a+1},K^{\{l\}})\left(\T_{m,l}^{(d)}\left(M'\cap(m..l),\frac1{\zeta-d-1}\right)\right)\biggr]\Biggr)\times\\[6pt]
\times\E_j(K^{(l,j)},L)\bigl(\T_{l,j}^{(d)}(M\cap(l..j),R)\bigl)\C(i,m)^{-1}
\pmod{\bar\I_{[i..j)}}.
\end{array}
$$

Applying~(\ref{equation:low5:8.5}), Proposition~\ref{proposition:low5:3}
and Lemma~\ref{lemma:low5:2}, we obtain that
the product in the first pair of the square brackets of the
above formula equals
$\E_j(K-\epsilon_{a+1},L)\bigl(\S_{i,m}^d\T^{(d)}_{m,j}(M',R)\bigl)$
modulo $\bar\I_{[i..j)}$ and the product in the second pair of
the square brackets equals $\E_l(K^{(i,l)}-\epsilon_{a+1},K^{\{l\}})
\Bigl(\S_{i,m}^d\T^{(d)}_{m,l}\Bigl(M'$ $\cap(m..l),\frac1{\zeta-d-1}\Bigr)\Bigr)$
modulo $\bar\I_{[i..l)}$. Now application of
Lemma~\ref{lemma:low5:1}\ref{lemma:low5:1:part:1},
Proposition~\ref{proposition:low5:2} and
Definition~\ref{definition:low5:2} concludes the proof.
\endproof

\section{Coefficients}\label{Coefficients}

\begin{proposition}\label{proposition:coeff:2}
Let $1\le i\le j\le n$, $c_1,\ldots,c_{n-1}$ be nonnegative integers
and $I$ denote the left ideal of $\U(n)$,
generated by $\E_t^{(c_t+1)}\U^+(n)$ for $t\in [i..j)$.
Assuming additionally $c_0=0$, we have
$$
\left(\prod_{t=1}^{n-1}\E_t^{(c_t+\delta_{t\in[i..j)})}\right)\S_{i,j}
{\equiv }\!\left(\prod_{t\in[i..j)}
(\B(i,t){+}c_{i-1}{-}c_i)
\right)\!
\E_1^{(c_1)}\cdots\E_{n-1}^{(c_{n-1})}\!\!\!\!\pmod I.\!
$$
\end{proposition}
\Proof This equivalence was actually proved in the course of the proof
of~\cite[Proposition 4.5]{Kleshchev_gjs11}. \endproof

\begin{proposition}\label{proposition:coeff:3}
Let $1\le i_1\le j_1\le n$, \ldots, $1\le i_m\le j_m\le n$ and
$c_1,\ldots,c_{n-1}$ be nonnegative integers.
We put $c_0:=0$,
$c^{(s)}_t:=c_t+\delta_{t\in[i_{s+1}..j_{s+1})}$\linebreak$+\cdots+\delta_{t\in[i_m..j_m)}$
for
$t=0,\ldots,n-1$
and $C:=(c_1,\ldots,c_{n-1})$. We have
$$
\begin{array}{l}
\displaystyle \left(\prod_{t=1}^{n-1}\E_t^{(c^{(0)}_t)}\right)
\S_{i_1,j_1}\cdots\S_{i_m,j_m}\\[20pt]
\displaystyle\equiv
\left(\prod_{s=1}^m\prod_{t\in[i_s..j_s)}
(\B(i_s,t)               +c^{(s)}_{i_s-1}-c^{(s)}_{i_s})
\right)
\left(\prod_{t=1}^{n-1}\E_t^{(c_t)}\right)
\pmod{\J^{(C)}}.
\end{array}
$$
\end{proposition}
\Proof We apply induction on $m$. Since the case $m=0$ is obvious,
we assume that $m>0$. The inductive hypothesis implies
$$
\begin{array}{l}
\displaystyle \left(\prod_{t=1}^{n-1}\E_t^{(c^{(1)}_t)}\right)
\S_{i_2,j_2}\cdots\S_{i_m,j_m}\\[20pt]
\displaystyle\equiv
\left(\prod_{s=2}^m\prod_{t\in[i_s..j_s)}
(\B(i_s,t)+c^{(s)}_{i_s-1}-c^{(s)}_{i_s})
\right)
\left(\prod_{t=1}^{n-1}\E_t^{(c_t)}\right)
\pmod{\J^{(C)}}.
\end{array}
$$

On the other hand, Proposition~\ref{proposition:coeff:2} implies
$$
\left(\prod_{t=1}^{n-1}\E_t^{(c^{(0)}_t)}\right)\S_{i_1,j_1}
{\equiv }\left(\prod_{t\in[i_1..j_1)}
(\B(i_1,t){+}c^{(1)}_{i_1-1}{-}c^{(1)}_{i_1})\right)
\left(\prod_{t=1}^{n-1}\E_t^{(c^{(1)}_t)}\right)\!\!\!\!\!\!
\pmod I,
$$
where $I$ is the left ideal of $\U(n)$ generated by
$\E_t^{(c^{(0)}_t)}\U^+(n)$ for $t\in [i_1..j_1)$.
Finally, it remains to notice the
inclusion
$I\S_{i_2,j_2}\cdots\S_{i_m,j_m}\subset\J^{(C)}$. \endproof

For a sequence $C=(c_1,\ldots,c_{n-1})$ of integers and integers
$k,i,t$ such that $1\le i,t<n$, we put for brevity
$\B^{C,k}(i,t):=\B(i,t)+c_{i-1}-c_i+\delta_{t\ge k}(c_{t+1}-c_t)$,\label{BCk}
where $c_0=c_n=0$.

To formulate the next lemma, we need to introduce rational expressions
similar to $\xi_{r,s}(M)$ of~\cite[Definition 2.9]{Kleshchev2}.
For any integers $i$ and $j$, sequences of integers
$C=(c_1,\ldots,c_{n-1})$, $K=(k_1,\ldots,k_d)$
$L=(l_1,\ldots,l_q)$, a subset $M\subset(i..j)$
and a rational expression $R\in\Z(\zeta)$ such that
$1\le i<j\le n$, $i\le k_1\le\cdots\le k_d\le j$ and $q=d$ or $q=d+1$,
we define $\rho^{(C)}(i,j,K,L,M,R)\in\bar\U^0(n)$
inductively on $|M|$ as follows,
additionally assuming $c_0=0$ and $c_n=0$.

{\it Case $M=\emptyset$.} We put
\begin{equation}\label{equation:coeff:2.5}
\begin{array}{l}
\displaystyle\rho^{(C)}(i,j,K,L,\emptyset,R){:=}
\displaystyle\prod_{s=1}^d\Biggl(\prod_{t\in[i..j)}(d{-}s{+}1)^{\delta_{t=j-1\ge k_s}}\Bigl(\B^{C,k_s}(i,t){-}d\\
+s{+}\delta_{t=j-1\ge k_s}\cdot\sum_{s<h\le q}l_h\Bigr)\Biggr)
R\Bigl(\C(i,j){+}c_{i-1}{-}c_i{-}c_{j-1}{+}c_j{+}\smallsum L\Bigr).
\end{array}
\end{equation}

{\it Case $M\ne\emptyset$.}
Let $m:=\min M$ and $M':=M\setminus\{m\}$. We put
\begin{equation}\label{equation:coeff:2.75}
\begin{array}{l}
\displaystyle\rho^{(C)}(i,j,K,L,M,R):=\Biggl(\rho^{(C)}(i,j,K,L,M',R)-
\zeta^{(C)}(i,m,K)\times\\[12pt]
\times\rho^{(C)}(m,j,K^{(m,j)},L,M',R)\Biggr)
\Bigl(\C(i,m){+}c_{i-1}{-}c_i{-}c_{m-1}{+}c_m\Bigr)^{-1},
\end{array}
\end{equation}
where
$
\displaystyle\zeta^{(C)}(i,m,K):=\prod_{s=1}^d\prod_{t\in[i..m)}(d-s+1)^{\delta_{t=m-1\ge k_s}}\times $\\
$\displaystyle\times\Bigl(
\B^{C,k_s}(i,t)+(\delta_{t<m-1}+\delta_{t=m-1<k_s})(s-d)\Bigr)$.\\

The reason for introducing these expresions is explained
as follows.

\begin{lemma}\label{lemma:coeff:1}
Let $1\le i<j\le n$, $d\ge1$, $M\subset(i..j)$ and $R\in\Z(\zeta)$.
Let $K=(k_1,\ldots,k_d)$, $L=(l_1,\ldots,l_q)$ and $C=(c_1,\ldots,c_{n-1})$
be sequences of integers such that all entries of $C$ are nonnegative,
$i\le k_1\le\cdots\le k_d\le j$ and $q=d$ or $q=d+1$.
We put $\bar c_t:=c_t+\sum_{s=1}^d\delta_{t\in[i..k_s)}$.
Modulo $\bar\J^{(C)}+\bar\I_{[\min K..j)}$, we have
$$
\left(\prod_{t=1}^{n-1}\E_t^{(\bar c_t)}\right)
\E_j(K,L)\left(\T_{i,j}^{(d)}(M,R)\right)\equiv
\rho^{(C)}(i,j,K,L,M,R)
\left(\prod_{t=1}^{n-1}\E_t^{(c_t)}\right).
$$
\end{lemma}
\Proof
We apply induction on $|M|$.
By Lemma~\ref{lemma:low5:2.5}, we have
$$
\begin{array}{l}
\displaystyle\E_j(K,L)\bigl(\T_{i,j}^{(d)}(\emptyset,R)\bigl)\equiv \S_{i,K}d^{\underline r}\times\\[6pt]
\displaystyle \times
\prod_{s=1}^r\left(
\left(\B(i,j-1)-d+s+{\textstyle\sum_{s<h\le q}}l_h\right)
\prod_{t\in[k_s..j-1)}\B(i,t)-d+s
\right)\times\\[20pt]
\times R(\C(i,j)+\smallsum L)\pmod{\bar\I_{[\min K..j)}},
\end{array}
$$
where $r$ is the number of entries of $K$ less than $j$.
Now Proposition~\ref{proposition:coeff:3} implies
\begin{equation}\label{equation:coeff:3}
\begin{array}{l}
\displaystyle\left(\prod_{t=1}^{n-1}\E_t^{(\bar c_t)}\right)\S_{i,K}\equiv
\prod_{s=1}^d\left[\prod_{t\in[i..k_s)}\left(\B(i,t)+c_{i-1}-(c_i+q_s)\right)\right]\times\\[20pt]
\displaystyle\times\left(\prod_{t=1}^{n-1}\E_t^{(c_t)}\right)
\pmod{\J^{(C)}},
\end{array}
\end{equation}
where $q_s$ is the number of elements of the sequence $k_{s+1},\ldots,k_d$
that are greater than $i$.
Let $a$ denote the number of entries of $K$
equal to $i$.
If $s\le a$ then the product in the square brackets
of~(\ref{equation:coeff:3}) is empty.
However, if $s>a$ then $q_s=d-s$.
Thus~(\ref{equation:coeff:3}) allows the following reformulation
$$
\begin{array}{l}
\displaystyle\left(\prod_{t=1}^{n-1}\E_t^{(\bar c_t)}\right)\S_{i,K}\equiv
\left(\prod_{s=1}^d\prod_{t\in[i..k_s)}(\B(i,t)+c_{i-1}-c_i-d+s)\right)
\times\\[20pt]
\times\displaystyle
\left(\prod_{t=1}^{n-1}\E_t^{(c_t)}\right)
\pmod{\J^{(C)}},
\end{array}
$$
Now throwing the corresponding elements of $\bar\U^0(n)$
over $\left(\prod_{t=1}^{n-1}\E_t^{(c_t)}\right)$
and applying~(\ref{equation:coeff:2.5}), we obtain the required formula.

Now suppose that $M\ne\emptyset$ and that the lemma is true
for sets of smaller cardinality.
We put $m:=\min M$ and $M':=M\setminus\{m\}$.
In the present case, equivalence~(\ref{equation:low5:2.5}) also holds,
where $P=\bar\I_{[\min K..j)}$. The inductive hypothesis implies
$$
\begin{array}{l}
\displaystyle
\left(\prod_{t=1}^{n-1}\E_t^{(\bar c_t)}\right)
\E_j(K,L)\bigl(
\T_{i,j}^{(d)}(M',R)
\bigr)\equiv \\[12pt]
\displaystyle
\rho^{(C)}(i,j,K,L,M',R)
\left(\prod_{t=1}^{n-1}\E_t^{(c_t)}\right)
\pmod{\bar\J^{(C)}+\bar\I_{[\min K..j)}}.
\end{array}
$$
We put $\hat c_t:=c_t+\sum_{s=1}^d\delta_{t\in[m..k_s)}$ and
$\hat C:=(\hat c_1,\ldots,\hat c_{n-1})$.
Clearly, $\bar c_t=\hat c_t+\sum_{s=1}^d\delta_{t\in[i..k_s^{(i,m)})}$.
Now the inductive hypothesis implies
$$
\begin{array}{l}
\displaystyle
\left(\prod_{t=1}^{n-1}\E_t^{(\bar c_t)}\right)
\E_m(K^{(i,m)},K^{\{m\}})\bigl(\S_{i,m}^d\bigl)
\equiv \\[12pt]
\displaystyle
\rho^{(\hat C)}(i,m,K^{(i,m)},K^{\{m\}},\emptyset,1)
\left(\prod_{t=1}^{n-1}\E_t^{(\hat c_t)}\right)
\pmod{\bar\J^{(\hat C)}+\bar\I_{[\min K..m)}}.
\end{array}
$$
One can easily observe that
$
\rho^{(\hat C)}(i,m,K^{(i,m)},K^{\{m\}},\emptyset,1)=
\zeta^{(C)}(i,m,K).
$
(the last element was introduced exactely for this equality).
Since by the inductive hypothesis we have
$$
\begin{array}{l}
\displaystyle\left(\prod_{t=1}^{n-1}\E_t^{(\hat c_t)}\right)
\E_j(K^{(m,j)},L)\bigl(\T_{m,j}^{(d)}(M',R)\bigl)\equiv \\[12pt]
\displaystyle\rho^{(C)}(m,j,K^{(m,j)},L,M',R)
\left(\prod_{t=1}^{n-1}\E_t^{(c_t)}\right)
\pmod{\bar\J^{(C)}+\bar\I_{[\min K..j)}},
\end{array}
$$
it remains to notice the following quite obvious formulae

\begin{tabular}{llr}
\hspace*{44pt}&$\bar\I_{[\min K..m)}\E_j(K^{(m,j)},L)\bigl(\T_{m,j}^{(d)}(M',R)\bigl)\in\bar\I_{[\min K..j)}$ &\\[6pt]
&             $\bar\J^{(\hat C)}\E_j(K^{(m,j)},L)\bigl(\T_{m,j}^{(d)}(M',R)\bigl)\in\bar\J^{(C)}$.&\hspace{38pt} $\square$
\end{tabular}

\begin{lemma}\label{lemma:coeff:2}
Let $\X\in\U^{-,0}(n)$
be an element of weight
$\sigma=-m_1\alpha_1-\cdots-m_{n-1}\alpha_{n-1}$ and
$C=(c_1,\ldots,c_{n-1})$ be a sequence of nonnegative integers.
Tnen there exists some $\rho\in\U^0(n)$ such that
$$
\left(\prod_{t=1}^{n-1}\E_t^{(c_t+m_t)}\right)\X
\equiv \rho \left(\prod_{t=1}^{n-1}\E_t^{(c_t)}\right)
\pmod{\J^{(C)}}.
$$
\end{lemma}
\Proof
Take any integers $t=1,\ldots,n-1$ and $i,j,r,c$ such that $1\le i<j\le n$
and $r,c\ge0$.
Using~(\ref{equation:algebras:0.5}), we see that
$\E_t^{(c)}\F_{i,j}^{(r)}=\sum_{s\ge0}F_s\E_t^{(c-s)}$,
where $F_s$ is an element of $\U^{-,0}(n)$ having weight
$-r\alpha(i,j)+s\alpha_t$.

Now notice that it suffices to prove the lemma for $\X=\F^{(N)}$,
where $N\in UT(n)$. Applying the formula we have just proved, we obtain
$$
\left(\prod_{t=1}^{n-1}\E_t^{(c_t+m_t)}\right)\F^{(N)}
=\sum_{s_1,\ldots,s_{n-1}\ge0}F_{s_1,\ldots,s_{n-1}}
\left(
\prod_{t=1}^{n-1}\E_t^{(c_t+m_t-s_t)}
\right),
$$
where $F_{s_1,\ldots,s_{n-1}}$ is an element of
$\U^{-,0}(n)$ having weight $\sigma+s_1\alpha_1+\cdots+s_{n-1}\alpha_{n-1}$.
In particular, $F_{s_1,\ldots,s_{n-1}}=0$ if
some $s_t>m_t$ and
$F_{m_1,\ldots,m_{n-1}}$ belongs to $\U^0(n)$, since this element
has weight zero. The proof concludes with noticing that
$
\prod_{t=1}^{n-1}\E_t^{(c_t+m_t-s_t)}\in\J^{(C)}
$
if some $s_t<m_t$.
\endproof

The next lemma follows directly from the definition of
$\rho^{(C)}(i,j,K,L,M,R)$.

\begin{lemma}\label{lemma:coeff:3}
If the entries at positions $i-1,\ldots,j$ of $C$ and $C'$ coinside,
then $\rho^{(C)}(i,j,K,L,M,R)=\rho^{(C')}(i,j,K,L,M,R)$.
Let $a\in\Z$, $j<n$, $L'=L\conc(a)$ if $|L|=d$ and
$L'=L+a\epsilon_{d+1}$ if $|L|=d+1$.
We have
$$
\rho^{(C+a\epsilon_j)}(i,j,K,L,M,R)=\rho^{(C)}(i,j,K,L',M,R).
$$
\end{lemma}

\section{Integral elements}\label{Integral elements}

In this section, for $N\in UT(n)$ we denote by $N_t$ the sum
$\sum_{1\le a\le n}N_{a,t}$ (the sum of elements in column $t$ of $N$).

\begin{proposition}[\mbox{\cite[2.9]{Carter_Lusztig}}]\label{proposition:int:1}
Let $1\le i\le j\le n$ and $d\ge 1$. Then we have
$$
\S_{i,j}^d=\sum_N\F^{(N)}\prod_{i<t\le j}N_t!\C(i,t)^{\underline{d-N_t}},
$$
where the summation runs over all $N\in UT(n)$ such that
$\F^{(N)}$ has weight $-d\alpha(i,j)$.
\end{proposition}
In other words, the summation runs over all $N\in UT(n)$
such that\linebreak $\sum_{1\le a\le t<b\le n}N_{a,b}=d\delta_{i\le t<j}$
for any $t=1,\ldots,n-1$. Clearly, if $i<j$ then $N_j=d$ for any such $N$.

Let $X\in\bar\U^{-,0}(n)$. By Lemma~\ref{lemma:algebras:1},
we have a unique representation
$
X=\sum_{N\in UT(n)}\F^{(N)}H_N,
$
where $H_N\in\bar\U^0(n)$. In that case, $H_N$ is called the
\linebreak
{\it $\F^{(N)}$-coefficient} of X. The proof of the following lemma
is similar to that of~\cite[Lemma~2.4]{Kleshchev2}.

\begin{lemma}\label{lemma:int:1}
For any integers $1\le i<j\le n$, $d\ge1$ and set $M\subset(i..j)$,
we have
$\T_{i,j}^{(d)}(M,1)\in\U^-(i,j)\mathcal U^0(i,j-1)$.
\end{lemma}
\Proof We fix $d\ge1$.
For any matrix $N\in UT(n)$ and any nonempty set $M$ such that
$\F^{(N)}$ has weight $-d\alpha(i,j)$ and $M\subset(i..j)$,
where $i$ and $j$ are some (uniquely determined) integers,
we define the polynomial $P_{N,M}\in\mathcal U^0(\min M,j-1)[x]$.
We require the two conditions
{
\renewcommand{\labelenumi}{{\rm \theenumi}}
\renewcommand{\theenumi}{{\rm(\arabic{enumi})}}
\begin{enumerate}
\item\label{lemma:int:1:condition:1}
      $P_{N,M}=P_{K,M}$ if $N_{a,b}=K_{a,b}$ for all $1\le a\le n$ and
      $\min M\le b\le n$;
\item\label{lemma:int:1:condition:2}
      for any $1\le i<j\le n$, nonempty $M\subset(i..j)$ and
      $N\in UT(n)$ such that $\F^{(N)}$ has weight $-d\alpha(i,j)$,
      the $\F^{(N)}$-coefficient of $\T_{i,j}^{(d)}(M,1)$ is
      $$
        d!\left(\prod_{i<t<\min M}N_t!\C(i,t)^{\underline{d-N_t}}\right)P_{N,M}(\C(i,\min M)).
      $$
\end{enumerate}
}
In view of Proposition~\ref{proposition:int:1}, constructing such polynomials
will automatically prove the lemma. Indeed, if $M=\emptyset$ then
$\T_{i,j}^{(d)}(M,1)=\S_{i,j}^d\in$\linebreak
$\U^-(i,j)\mathcal U^0(i,j-1)$ by Definition~\ref{definition:low5:2} and
Proposition~\ref{proposition:int:1}. If $M\ne0$ then by
condition~\ref{lemma:int:1:condition:2}, the $\F^{(N)}$-coefficient
of $\T_{i,j}^{(d)}(M,1)$ is an integer polynomial in $\H_i,\ldots,\H_{j-1}$,
that is an element of $\mathcal U^0(i,j-1)$.

We apply induction on $|M|>0$. Put $m:=\min M$, $M':=M\setminus\{m\}$ and
$m':=\min M'$ if $M'$ is not empty.
We take a matrix $N\in UT(n)$ such that $\F^{(N)}$ has weight
$-d\alpha(i,j)$ for some $i$ and $j$ with $M\subset(i..j)$.
We are going to define $P_{N,M}$ so that
condition~\ref{lemma:int:1:condition:2} holds.

{\it Case $N_m<d$}.
Since the $\F^{(N)}$-coefficient of
$\S_{i,m}^d\T_{m,j}^{(d)}(M',1)$ equals zero in this case,
we can define
$$
P_{N,M}:=N_m!x^{\underline{d-N_m-1}}
\left(\prod_{m<t<m'}N_t!(x+\C(m,t))^{\underline{d-N_t}}\right)
P_{N,M'}(x+\C(m,m'))
$$
if $M'\ne\emptyset$ and
$$
P_{N,M}:=N_m!x^{\underline{d-N_m-1}}
\left(\prod_{m<t<j}N_t!(x+\C(m,t))^{\underline{d-N_t}}\right)
$$
if $M'=\emptyset$.

{\it Case $N_m=d$.}
We have the decomposition $N=N'+N''$,
where $N',N''\in UT(n)$, $N'_{a,b}=\delta_{b\le m}N_{a,b}$ and
$N''_{a,b}=\delta_{a\ge m}N_{a,b}$.
This decomposition follows from $\sum_{a\le m-1<b}N_{a,b}=d$
and the fact that all entries of $N$ are nonnegative.
Thus we have $\F^{(N)}=\F^{(N')}\F^{(N'')}$.
If $M'\ne\emptyset$ then by condition
\ref{lemma:int:1:condition:1}
and the inductive hypothesis
we get $P_{N,M'}=P_{N'',M'}$.
Therefore, we can define
$$
\begin{array}{l}
\displaystyle P_{N,M}:=\frac{d!}x\left(\left(\prod_{m<t<m'}N_t!(x+\C(m,t))^{\underline{d-N_t}}\right)P_{N,M'}(x+\C(m,m'))\right.\\[12pt]
\displaystyle \left.-\left(\prod_{m<t<m'}N_t!\C(m,t)^{\underline{d-N_t}}\right)P_{N,M'}(\C(m,m'))\right)
\end{array}
$$
if $M'\ne\emptyset$ and
$$
\displaystyle P_{N,M}:=\frac{d!}x
\left(\left(\prod_{m<t<j}N_t!(x+\C(m,t))^{\underline{d-N_t}}\right)
-\left(\prod_{m<t<j}N_t!\C(m,t)^{\underline{d-N_t}}\right)\right)
$$
if $M'=\emptyset$.
One can easily see that $P_{N,M}\in\mathcal U^0(\min M,j-1)[x]$.
Moreover, condition~\ref{lemma:int:1:condition:1} inductively follows
from the formulae defining $P_{N,M}$.
\endproof

Each commutative ring $\mathcal U^0(a,b)$ is generated by elements
$\H_a,\ldots,\H_b$ freely over $\Z$.
Therefore, $\mathcal U^0(a,b)$ is isomorphic to $\Z[x_a,\ldots,x_b]$
and is a UFD. This fact is used below.

\begin{corollary}\label{corollary:int:1}
In Lemmas~\ref{lemma:low5:3} and~\ref{lemma:low5:4},
the expressions in
square brackets
belong to
$\mathcal U^{-,0}(i,l)\oplus\bar\I_{[\min K..l)\cup\{l-1\}}$,
whenever the corresponding case occurs.
\end{corollary}
\Proof Denote this expression by $\bar X$ and
let $X$ be the element of $\bar\U^{-,0}(i,l)$ such that
$\bar X\equiv X\pmod{\bar\I_{[\min K..l)\cup\{l-1\}}}$.
Without loss of generality
we can assume that $R=1$ and $M\cap(l..j)=\emptyset$.
Let $\bar Y:=\E_j(K^{(l,j)},L)\bigl(\T_{l,j}^{(d)}(\emptyset,1)\bigl)$.
Thus $\bar Y$ is the second factor in the product where $\bar X$ occurs.
Let $Y$ be the element of $\bar\U^{-,0}(l,j)$ such that
$\bar Y\equiv Y\pmod{\I_{[l..j)}}$.
By Lemma~\ref{lemma:coeff:1}, we have $Y\ne0$,
since $\rho^{(0^{n-1})}(l,j,K^{(l,j)},L,\emptyset,1)\ne0$.
On the other hand, $Y\in\mathcal U^{0,-}(l,j)$.
By Lemma~\ref{lemma:int:1}, applying Lemma~\ref{lemma:low5:3} or
Lemma~\ref{lemma:low5:4} respectively, we get $XY\in\mathcal U^{0,-}(i,j)$.
Let $D(i,l)$ denote the set of all products (including empty) of elements of
the form $\C(s,t)+N$, where $i\le s<t\le j$ and $N\in\Z$.

Choose some $N\in UT(n)$ such that the $\F^{(N)}$-coefficient $H_N$
of $Y$ is not equal to zero.
Let $M\in UT(n)$ be a matrix such that the $\F^{(M)}$-coefficient $H_M$
of $X$ is not equal to zero.
We have $\F^{(M)}\F^{(N)}=\F^{(M+N)}$. Thus the $\F^{(M+N)}$-coefficient
of $XY$ is $\tau_N(H_M)H_N$ by the first formula of~(\ref{equation:algebras:1})
and the remark before Lemma~\ref{lemma:algebras:1}.
We have
$$
H_N\in\mathcal U^0(l,j),\quad \tau_N(H_M)H_N\in\mathcal U^0(i,j),\quad
H_M=h/f
$$
for some $h\in\mathcal U^0(i,l)$ and $f\in D(i,l)$.
The above representation for $H_M$ is derived from
Definition~\ref{definition:low5:2} and
Proposition~\ref{proposition:low5:-1&low5:2.5}\ref{proposition:low5:-1}.
We have $\tau_N(H_M)=\tau_N(h)/\tau_N(f)$,
$\tau_N(h)\in\mathcal U^0(i,l)$ and $\tau_N(f)\in D(i,l)$.
Therefore, $\tau_N(f)$ divides $\tau_N(h)H_N$ in the ring $\mathcal U^0(i,j)$.
Since $H_N\in\mathcal U^0(l,j)$, the prime decomposition of $H_N$
can not contain factors of the form $\C(s,t)+N$, where $i\le s<t\le l$.
Therefore, $\tau_N(H_M)=\tau_N(h)/\tau_N(f)\in\mathcal U^0(i,l)$,
whence $H_M\in\mathcal U^0(i,l)$. To prove the latter,
one needs to apply the inverse map to $\tau_N$.
\endproof

In the next lemma, we
use the following notation:
$C=(c_1,\ldots,c_{n-1})$, $c_0=c_n=0$,
$K=(k_1,\ldots,k_d)$, $L=(l_1,\ldots,l_q)$,
where $d\ge1$ and $q=d$ or $q=d+1$.

\begin{lemma}\label{lemma:coeff:0.5}
We have $\rho^{(C)}(i,j,K,L,M,R)\in\mathcal U^0(i,j)$ if $R=1$ or
if $R=1/(\zeta-d)$, $K\ne(j^d)$ and $l_1=0$.
\end{lemma}
\Proof We shall only sketch the proof,
leaving technical details to the reader.
For $m\in(i..j)$, let $\sigma_m$ be the endomorphism of the
ring $\mathcal U^0(i,j)$ defined by
$\sigma_m(\H_i):=i-m+\H_m-c_{i-1}+c_i+c_{m-1}-c_m$ and
$\sigma_m(\H_t):=\H_t$ for $t\ne i$.
We shall prove inductively on $|M|$ the two facts
\begin{enumerate}
\item\label{lemma:coeff:0.5:part:1} $\rho^{(C)}(i,j,K,L,M,R)\in\mathcal U^0(i,j)$;
\item\label{lemma:coeff:0.5:part:2} $\sigma_m\Bigl(\rho^{(C)}(i,j,K,L,M,R)-
      \zeta^{(C)}(i,m,K)\rho^{(C)}(m,j,K^{(m,j)},L,M,R)\Bigr)=0$
      for any $m\in(i..\min M\cup\{j\})$.
\end{enumerate}

The case $M=\emptyset$ follows from the definition of
$\rho^{(C)}(i,j,K,L,\emptyset,R)$ and direct calculations.

Now let $M\ne\emptyset$,
$m'=\min M$ and $M'=M\setminus\{m'\}$.
To prove condition~\ref{lemma:coeff:0.5:part:1},
we first note that by the inductive hypothesis the element
\begin{equation}\label{equation:coeff:2.875}
\rho^{(C)}(i,j,K,L,M',R)-\zeta^{(C)}(i,m,K)\rho^{(C)}(m,j,K^{(m,j)},L,M',R)
\end{equation}
belongs to $\mathcal U^0(i,j)$. Thus it can be considered as
a polynomial in $\H_i$ over $\mathcal U^0(i+1,j)$.
Applying $\sigma_{m'}$ to~(\ref{equation:coeff:2.875}) and
using the inductive hypothesis, we obtain that
$i-m'+\H_{m'}-c_{i-1}+c_i+c_{m'-1}-c_{m'}$ is a root of this polynomial.
Therefore, it is divisible by $\C(i,m')+c_{i-1}-c_i-c_{m'-1}+c_{m'}$ and
the coefficients of the quotient belong to $\mathcal U^0(i+1,j)$.
The required result follows now from the definition of $\rho^{(C)}(i,j,K,L,M,R)$.
Condition~\ref{lemma:coeff:0.5:part:2} can be checked by direct calculations.
\endproof

\begin{proposition}\label{proposition:commpoly:1}
Let $f_1,\ldots,f_a$ be first degree polynomials of $\mathcal U^0(i,j)$
having the form $f_h=\H_{m_h}+g_h$, where $g_h$ is a $\Z$-linear combination
of the unit and the variables $\H_t$ for $t>m_h$ and
$i\le m_1<\cdots<m_a\le j$. Let $I$ denote the ideal of $\mathcal U^0(i,j)$
generated by $f_1,\ldots,f_a$.
Then
\begin{enumerate}
\item\label{proposition:commpoly:1:case:1}
 $I$ is a prime ideal;
\item\label{proposition:commpoly:1:case:2}
 a first degree polynomial belongs to $I$ if and only if it is
 a $\Z$-linear combination of $f_1,\ldots,f_a$.
\end{enumerate}
\end{proposition}

The next result states the property of
$\rho^{(C)}(i,j,K,(0^{d+1}),M,1)$ similar to that of $\xi_{r,s}(M)$
proved \mbox{in~\cite[Lemma 2.11]{Kleshchev2}}.

\begin{lemma}\label{lemma:commpoly:1}
Let $1\le i<j\le n$, $d\ge1$, $C=(c_1,\ldots,c_{n-1})$ and
$K=(k_1,\ldots,k_d)$ be sequences of integetrs and $M=\{m_1,\ldots,m_a\}$
be a set such that $i\le k_1\le\cdots\le k_d\le j$ and $i<m_1<\cdots<m_a<j$.
Suppose that $h\mapsto(t_h,s_h)$ is an injective map from
$[1..a]$ to $[i..j)\times[1..d]$ such that $t_h\ge m_h$ for all $h$.
Modulo the ideal of $\mathcal U^0(i,j)$ generated by
$
\displaystyle\B^{C,k_{s_h}}(m_h,t_h)-d+s_h,
$
where $h=1,\ldots,a$, we have
$$
\begin{array}{l}
\displaystyle\rho^{(C)}(i,j,K,(0^{d+1}),M,1)\equiv \\[12pt]
\displaystyle d^{\underline r}\prod\Bigl\{
\B^{C,k_s}(i,t)-d+s:
(t,s)\in[i..j)\times[1..d]\setminus\{(t_h,s_h):h=1,\ldots,a\}\Bigl\},
\end{array}
\!\!\!
$$
where $r$ is the number of entries of $K$ less than $j$.
\end{lemma}
\Proof
We apply induction on $|M|$.
The case $M=\emptyset$ follows immediately from the definition.
Now let $M\ne\emptyset$.
We put $M':=\{m_2,\ldots,m_a\}$, $X:=[i..j)\times[1..d]$,
$Y:=[m_1..j)\times[1..d]$, $\Phi:=\{(t_h,s_h):h=1,\ldots,a\}$ and
$\Psi:=\{(t_h,s_h):h=2,\ldots,a\}$. The inductive hypothesis implies
$$
\begin{array}{l}
\arraycolsep=0pt
\rho^{(C)}(i,j,K,(0^{d+1}),M',1)=
\displaystyle d^{\underline r}
\left(\prod\nolimits_{(t,s)\in X\setminus\Psi}
\B^{C,k_s}(i,t)-d+s\right)
+f,\\[12pt]
\displaystyle \rho^{(C)}(m_1,j,K^{(m_1,j)},(0^{d+1}),M',1){=}d^{\underline r}
\left(\prod\nolimits_{(t,s)\in Y\setminus\Psi}
\B^{C,k_s}(m_1,t){-}d{+}s\right){+}g,
\end{array}\!\!\!
$$
where $f$ and $g$ belong to the ideal $I$ of $\mathcal U^0(i,j)$
generated by
$
\B^{C,k_{s_h}}(m_h,t_h)-d+s_h
$
for $h=2,{\ldots},a$.
We put for brevity $x{:=}\C(i,m_1){+}c_{i-1}{-}c_i{-}c_{m_1-1}{+}c_{m_1}$.
Since $\B(i,t)+c_{i-1}-c_i=x+\B(m_1,t)+c_{m_1-1}-c_{m_1}$, we get
\begin{equation}\label{equation:commpoly:1}
\!\!\!\!
\!\!\!\!
\begin{array}{l}
\arraycolsep=0pt
\displaystyle x\rho^{(C)}(i,j,K,(0^{d+1}),M,1)=xd^{\underline r}
\left(\prod\nolimits_{(t,s)\in X\setminus\Phi}
(\B^{C,k_s}(i,t){-}d{+}s)\right){+}
\\[12pt]
\displaystyle
d^{\underline r}(\B^{C,k_{s_1}}(m_1,t_1){-}d{+}s_1)
\left(\prod\nolimits_{s=1}^d\prod\nolimits_{t\in[i..m_1-1)}
(\B^{C,k_s}(i,t){-}d{+}s)\right){\times}\\[12pt]
\displaystyle\times\Biggl[\left(\prod\nolimits_{s=1}^d
(x-1+(1-\delta_{m_1-1\ge k_s})(c_{m_1-1}-c_{m_1})-d+s)\right)
\times\\[12pt]
\displaystyle{\times}
\left(\prod\nolimits_{(t,s)\in Y\setminus\Phi}
(x{+}\B^{C,k_s}(m_1,t){-}d{+}s)\right){-}
\displaystyle\biggl(\prod\nolimits_{s=1}^d(d{-}s{+}1)^{\delta_{m_1-1\ge k_s}}
{\times}\!\!\!\!
\\[12pt]
\times
\displaystyle
\Bigl(x-1+(1-\delta_{m_1-1\ge k_s})(c_{m_1-1}-c_{m_1})
+\delta_{m_1-1<k_s}(s-d)\Bigl)\biggr)\times\\[12pt]
\displaystyle\times
\left(\prod\nolimits_{(t,s)\in Y\setminus\Phi}
(\B^{C,k_s}(m_1,t)-d+s)\right)\Biggl]+f+g',
\end{array}
\end{equation}
where $g'\in\mathcal U^0(i,j)g\subset I$.
The formal substitution $x\mapsto 0$ takes the expression
in the square brackets to zero, and therefore,
this expression is divisible by $x$.
Hence $f+g'=xq\in I$ for some $q\in\mathcal U^0(i,j)$.
By Proposition~\ref{proposition:commpoly:1}\ref{proposition:commpoly:1:case:1},
we have $x\in I$ or $q\in I$.
The former case is imposible by Proposition~\ref{proposition:commpoly:1}\ref{proposition:commpoly:1:case:2}.
Therefore $q\in I$. Now dividing~(\ref{equation:commpoly:1}) by $x$,
we prove the requred result.
\endproof

For $d\ge1$, $i<j$ and $M\subset(i..j)$, we define the following polymomials
$$\label{fijgij}
\begin{array}{l}
\arraycolsep=0pt
\displaystyle f_{i,j}^{(d)}(\emptyset){:=}\prod\nolimits_{t\in[i..j)}(y_{t+1}{-}x_i)^{\underline d},\,\,
g_{i,j}^{(d)}(\emptyset){:=}(y_j{-}x_i)^{\underline{d-1}}\prod\nolimits_{t\in[i..j-1)}(y_{t+1}{-}x_i)^{\underline d};\\[12pt]
f_{i,j}^{(d)}(M):=\tfrac{f_{i,j}^{(d)}(M')-f_{i,m}^{(d)}(\emptyset)f_{m,j}^{(d)}(M')}{x_m-x_i},\quad
g_{i,j}^{(d)}(M):=\tfrac{g_{i,j}^{(d)}(M')-f_{i,m}^{(d)}(\emptyset)g_{m,j}^{(d)}(M')}{x_m-x_i},\\[8pt]
\mbox{where $M\ne\emptyset$, $m=\min M$ and $M'=M\setminus\{m\}$.}
\end{array}
\!\!\!\!\!
$$
Similarly to Lemma~\ref{lemma:coeff:0.5}, one can prove that
$f_{i,j}^{(d)}(M)$ and $g_{i,j}^{(d)}(M)$ are elements of
$\Z[x_i,\ldots,x_{j-1},y_{i+1},\ldots,y_j]$.
For $M\subset(i..j)$, $l\in M$ and a strictly increasing sequence
$N=(i_1,\ldots,i_k)\subset M\cap(i..l)$, we put
$$\label{Gil}
G_{i,l}^{(d)}(M,N):=\prod_{r=0}^kg_{i_r,i_{r+1}}^{(d)}(M\cap(i_r..i_{r+1})),\;\;
\mbox{ where $i_0=i$ and $i_{k+1}=l$ }.
$$

\begin{lemma}\label{lemma:commpoly:2}
Let $M\subset(i..j)$ and $l\in M$. Then $f_{i,j}^{(d)}(M)$
belongs to the ideal of $\Z[x_i,\ldots,x_{j-1},y_{i+1},\ldots,y_j]$
generated by $G_{i,l}^{(d)}(M,N)$, where $N\subset M\cap(i..l)$.
\end{lemma}
\Proof
We apply induction on $|M|$.
For $M=\{l\}$, we have
$$
\begin{array}{l}
\displaystyle f_{i,j}^{(d)}(M):=
\displaystyle (y_l-x_i-d+1)g_{i,l}^{(d)}(\emptyset)\tfrac{\left(\prod\nolimits_{t\in[l..j)}(y_{t+1}-x_i)^{\underline d}\right)-\left(\prod\nolimits_{t\in[l..j)}(y_{t+1}-x_l)^{\underline d}\right)}{x_l-x_i},
\end{array}
$$
The substitution $x_i\mapsto x_l$ shows that the fraction in the right-hand side
of this expression is a polynomial of $\Z[x_i,\ldots,x_{j-1},y_{i+1},\ldots,y_j]$.

Now let $|M|>1$.
We put $m:=\min M$ and $M':=M\setminus\{m\}$.
First we consider the case $l=m$.
We put $m':=\min M'$.
By the inductive hypothesis, we have
$f_{i,j}^{(d)}(M')=g_{i,m'}^{(d)}(\emptyset)h$, where
$h\in\Z[x_i,\ldots,x_{j-1},y_{i+1},\ldots,y_j]$.
\linebreak
Hence $f_{i,j}^{(d)}(M')=g_{i,l}^{(d)}(\emptyset)h'$, where
$
h'=h
(y_l-x_i-d+1)(y_{m'}-x_i)^{\underline{d-1}}\times
$
\linebreak
$
\times\prod_{t\in[l..m'-1)}(y_{t+1}-x_i)^{\underline d}
$.
We have $f_{i,j}^{(d)}(M)=g_{i,l}^{(d)}(\emptyset)\bigl(h'-(y_l-x_i-d+1)f_{l,j}^{(d)}(M')\bigl)/(x_l-x_i)$.
Since $x_l-x_i$ and $g_{i,l}^{(d)}(\emptyset)$ are relatively prime, we have
$(h'-(y_l-x_i-d+1)f_{l,j}^{(d)}(M'))/(x_l-x_i)\in\Z[x_i,\ldots,x_{j-1},y_{i+1},\ldots,y_j]$.

Now we consider the case $m<l$.
The inductive hypothesis implies that $f_{i,j}^{(d)}(M')$ belongs to the ideal of
$\Z[x_i,\ldots,x_{j-1},y_{i+1},\ldots,y_j]$ generated by the polynomials
$G^{(d)}_{i,l}(M',N)$, where $N\subset M'\cap(i..l)$.
Since for $i_1\in M'$ one has
$$
\begin{array}{l}
g^{(d)}_{i,i_1}(M'\cap(i..i_1))=(x_m-x_i)g_{i,i_1}^{(d)}(M\cap(i..i_1))\\[6pt]
+(y_m-x_i-d+1)g_{i,m}^{(d)}(\emptyset)g_{m,i_1}^{(d)}(M\cap(m..i_1)),
\end{array}
$$
$f_{i,j}^{(d)}(M')$ belongs to the ideal $I$ of
$\Z[x_i,\ldots,x_{j-1},y_{i+1},\ldots,y_j]$,
generated by the polynomials
$(x_m-x_i)^{\delta_{m\notin N}}G_{i,l}^{(d)}(M,N)$, where $N\subset M\cap(i..l)$.
By the inductive hypothesis $f_{i,m}^{(d)}(\emptyset)f_{m,j}^{(d)}(M')\in I$.
Hence $f_{i,j}^{(d)}(M)(x_m-x_i)\in I$.
We have
$
f_{i,j}^{(d)}(M)(x_m-x_i)=\sum\nolimits_{N\subset M\cap(i..l)}(x_m-x_i)^{\delta_{m\notin N}}G^{(d)}_{i,l}(M,N)h_N
$
for some $h_N\in\Z[x_i,\ldots,x_{j-1},y_{i+1},\ldots,y_j]$.

For every subset $N\subset M\cap(i..l)$ such that $m\in N$, we consider the representation
$h_N=h'_N(x_m-x_i)+h''_N$,
where $h'_N\in\Z[x_i,\ldots,x_{j-1},y_{i+1},\ldots,y_j]$ and
$h''_N\in\Z[x_{i+1},\ldots,x_{j-1},y_{i+1},\ldots,y_j]$.
Hence the product
\begin{equation}\label{equation:commpoly:1.5}
g_{i,m}^{(d)}(\emptyset)\cdot\sum\nolimits_{N'\subset M'\cap(m..l)}G_{m,l}^{(d)}(M',N')h''_{N'}
\end{equation}
is divisible by $(x_m-x_i)$. Since $(x_m-x_i)$ does not divide the first factor
of~(\ref{equation:commpoly:1.5}), it divides the second factor, which however does not
depend on $x_i$. Therefore, the whole polynomial~(\ref{equation:commpoly:1.5}) equals zero.
Thus we have obtained
$$
{}\hspace{31pt}
f_{i,j}(M)=\sum_{\genfrac{}{}{0pt}{}{N\subset M\cap(i..l)}{m\notin N}}G^{(d)}_{i,l}(M,N)h_N+
\sum_{\genfrac{}{}{0pt}{}{N\subset M\cap(i..l)}{m\in N}}G_{i,l}^{(d)}(M,N)h'_N.
\hspace{31pt}\square
$$

\begin{corollary}\label{corollary:commpoly:1}
Let $C=(c_1,\ldots,c_{n-1})$ be a sequence of integers, $d\ge1$,
$1\le i<j\le n$, $M\subset(i..j)$ and $l\in M$.
Then $\rho^{(C)}(i,j,(j^d),(0^{d+1}),M,1)$ belongs to the ideal of
$\mathcal U^0(i,j)$ generated by
$$
\prod_{r=0}^k\tfrac1d\rho^{(C)}\left(i_r,i_{r+1},(i_{r+1}-1,i_{r+1}^{d-1}),(0^d,d),M\cap(i_r..i_{r+1}),\tfrac1{\zeta-d}\right),
$$
where $i=i_0<i_1<\cdots<i_k<i_{k+1}=l$, $k\ge0$ and $i_1,\ldots,i_k\in M$.
\end{corollary}
\Proof The result follows from the previous lemma and formulae~(\ref{equation:coeff:2.5})
and (~\ref{equation:coeff:2.75}) by the substitution
$y_t\mapsto t-1-\H_t$, $x_t\mapsto t-\H_t-c_{t-1}+c_t$. \endproof

\section{Proof of the main results}\label{Proof of the main results}

We define the hyperalgebra $U(n)$\label{U(n)}
over $\K$ to be $\U(n)\otimes_\Z\K$.
It is well known that every rational $\GL_n$-module
can be considered as a $U(n)$-module (see~\cite[I.7.11--I.7.16]{Jantzen1}).
We shall use the following notation:
$E_i^{(r)}:=\E_i^{(r)}\otimes 1_\K$,
$\binom{H_{i,j}}r:=\binom{\H_{i,j}}r\otimes 1_\K$,
$U^?(n):=\U^?(n)\otimes\K$, $U^?(i,j):=\U^?(i,j)\otimes\K$.
Any weight $\lambda=(\lambda_1,\ldots,\lambda_n)$ will also be identified
with the $\K$-algebra homomorphism $\pi_\lambda:U^0(n)\to \K$\label{pilm}
that takes $\binom{H_i}{r}$ to $\binom{\lambda_i}{r}+p\Z$.
A vector $v$ of a $U(n)$-module is called a {\it $U(n)$-high weight vector}
if it is a $U^{0}(n)$-weight vector and $E^{(r)}_iv=0$ for all $i=1,\ldots,n-1$ and $r>0$.
One can observe that a vector of a rational $\GL_n$-module
is a $\GL_n$-high weight vector if and only if it is
a $U(n)$-high weight vector.

For the rest of the paper, we fix $n>1$.
For $\lambda\in X^+(n)$, we denote by $\nabla_n(\lambda)$\label{CoWeyl}
the $\GL_n$-module
contravariantly dual to the Weyl module with highest weight $\lambda$.
It is well known that $\soc\nabla_n(\lambda)\cong L_n(\lambda)$.
We need the module $\nabla_n(\lambda)$ because we know all its
$U(n-1)$-high weight vectors.
Explicitly, for any $\mu\in X^+(n-1)$
such that
$\lambda_{i+1}\le\mu_i\le\lambda_i$ for $i=1,\ldots,n-1$,
there is a nonzero $U(n-1)$-high weight vector
$f_{\mu,\lambda}\in\nabla_n(\lambda)$\label{fmulm}
of weight $\mu$.
We shall express the above relation between $\mu$ and $\lambda$ by
$\mu\longleftarrow\lambda$.
These vectors have the property that
for any nonzero $U(n-1)$-high weight vector $v\in\nabla_n(\lambda)$
of weight $\mu\in X^+(n-1)$, there holds $\mu\longleftarrow\lambda$
and $v$ is a scalar multiple of $f_{\mu,\lambda}$
(see~\cite[Corollary 3.3]{Kleshchev_gjs11}).
We shall abbreviate $f_\lambda:=f_{\bar\lambda,\lambda}$\label{flm},
where $\bar\lambda=(\lambda_1,\ldots,\lambda_{n-1})$.
Thus $\soc\nabla_n(\lambda)$ is generated by $f_\lambda$ as a $\GL_n$-module.

\begin{proposition}\label{proposition:primvect:1}${}$
\!\!\!\!\!\!
We can choose the vectors $f_{\mu,\lambda}$ so that
$E_1^{(a_1)}{\cdots}E_{n-1}^{(a_{n-1})}f_{\mu,\lambda}$\linebreak $=f_\lambda$,
where
$a_i=\sum_{s=1}^i(\lambda_s-\mu_s)$.
\end{proposition}
\Proof
If all entries of $\lambda$ are nonnegative, then
this can be done by \cite[Lemma 2.6]{Kleshchev_gjs11}.
To this case, the arbitrary case can be reduced by
tensoring $\nabla_n(\lambda)$ with a sufficiently high power of the determinant
representation.
\endproof

Now we are going to associate with every homogeneous vector
$v\in\nabla_n(\lambda)$ the element $\c(v)\in\K$\label{cfv}
as follows.
Let $v$ have weight $\nu\in X(n-1)$ and
$a_i=\sum_{s=1}^i(\lambda_s-\nu_s)$. We choose $\c(v)$ so that
$E_1^{(a_1)}\cdots E_{n-1}^{(a_{n-1})}v=\c(v)f_\lambda$.
We define the relation $\dasharrow$ between vectors of $\nabla_n(\lambda)$
by the rules:
$v\dasharrow v$; $v\dasharrow E_l^{(r)}v$ for $1\le l<n-1$ and $r>0$;
if $v\dasharrow u$ and $u\dasharrow w$ then $v\dasharrow w$.

The next theorem is our main tool to find out whether
a vector of $\nabla_n(\lambda)$ is a nonzero $U(n-1)$-high weight vector.

\begin{theorem}\label{theorem:primvect:1}
Let a vector $v$ of $\nabla_n(\lambda)$ have weight $\nu\in X(n-1)$.
\begin{enumerate}
\item\label{theorem:primvect:1:part:1}
      $v=0$ if and only if $\c(u)=0$ for any $u$ such that $v\dasharrow u$.
\item\label{theorem:primvect:1:part:2}
      $v$ is a nonzero $U(n-1)$-high weight vector
      if and only if $\c(v)\ne0$ and $\c(u)=0$
      for any $u$ such that $u\ne v$ and $v\dasharrow u$.
\end{enumerate}
\end{theorem}
\Proof\ref{theorem:primvect:1:part:1}
Clearly, it suffices to prove that $v=0$ if $\c(u)=0$ for any $u$
such that $v\dasharrow u$. Assume this is wrong. Then there is
$\nu\in X(n-1)$ and a vector $v\in\nabla_n(\lambda)$ of weight $\nu$
such that $v\ne 0$ but $\c(u)=0$ for any $u$ such that $v\dasharrow u$.
We can assume that $\nu$ is the maximal weight with this property.

Take arbitrary $1\le l<n-1$ and $r>0$.
We have $v\dasharrow E_l^{(r)}v$.
If $E_l^{(r)}v\dasharrow u$ then $v\dasharrow u$ and thus $\c(u)=0$.
Since $E_l^{(r)}v$ has $U(n-1)$-weight strictly greater than $\nu$, we have
$E_l^{(r)}v$=0.

We have proved that $v$ is a nonzero $U(n-1)$-high weight vector
of $\nabla_n(\lambda)$.
By~\cite[Corollary 3.3]{Kleshchev_gjs11},
we have $\nu\longleftarrow\lambda$ and $v=\beta f_{\nu,\lambda}$
for some $\beta\in\K\setminus\{0\}$.
Multiplying this equality by $E_1^{(a_1)}\cdots E_{n-1}^{(a_{n-1})}$,
where $a_i=\sum_{s=1}^i(\lambda_s-\nu_s)$,
and taking into account Proposition~\ref{proposition:primvect:1},
we obtain $0=\c(v)f_\lambda=\beta f_\lambda$.
Hence $\beta=0$ contrary to assumption.

\ref{theorem:primvect:1:part:2}
Let $v$ be a nonzero $U(n-1)$-high weight vector.
For any vector $u$ such that $u\ne v$ and $v\dasharrow u$,
there exist $1\le l<n-1$ and $r>0$ such that $E_l^{(r)}v\dasharrow u$.
However, $E_l^{(r)}v=0$ and thus $u=0$. In particular, $\c(u)=0$.
If $\c(v)$ equaled zero, then by part~\ref{theorem:primvect:1:part:1}
we would get $v=0$, which is wrong.

Now let $\c(v)\ne0$ and $\c(u)=0$ for any $u$ such that
$u\ne v$ and $v\dasharrow u$. We have $v\ne0$.
Take any integers $l$ and $r$ such that $1\le l<n-1$ and $r>0$.
For any vector $u$ such that $E_l^{(r)}v\dasharrow u$,
we have $v\dasharrow u$ and $u\ne v$, since either $u=0$ or
$u$ has $U(n-1)$-weight strictly greater than $\nu$.
Thus $\c(u)=0$.
By part~\ref{theorem:primvect:1:part:1} we have $E_l^{(r)}v=0$.
\endproof

Denote by $S_{i,j}$,\label{SinU}
          $T_{i,j}^{(d)}(M,R)$,\label{TinU}
          $\xi^{(C)}(i,j,K,L,M,R)$\label{xi}
the images of $\S_{i,j}$,\linebreak $\T_{i,j}^{(d)}(M,R)$, $\rho^{(C)}(i,j,K,L,M,R)$
in $U(n)$ respectively, whenever the corresponding elements of $\bar\U(n)$
belong to $\U(n)$.
Let $I_S$ and $J^{(C)}$\label{IJinU}
be the images of the ideals $\I_S$ and $\J^{(C)}$
in $U(n)$ respectively.
The proof of the following result is similar to that of
Lemma~\ref{lemma:coeff:2}.

\begin{proposition}\label{lemma:commpoly:3}
Let $X$ be an element of $U^{-,0}(n)$ of $U(n)$-weight $\sigma$,
$1\le t<n$ and $r\ge0$. Then we have
$E_t^{(r)}X=\sum_{s=0}^bX_sE_t^{(r-s)}$,
where $X_s$ is an element of $U^{-,0}(n)$ of weight $\sigma+s\alpha_t$
and $b$ is the $\alpha_t$-coefficient of $-\sigma$.
\end{proposition}

In the sequel, we shall use the notation
$$\label{Bmulm}
B^{\mu,\lambda}(i,t):=t-i+\mu_i-\lambda_{t+1},\;\;\;
B^{\mu,\lambda,k}(i,t):=\left\{
\begin{array}{l}
t-i+\mu_i-\mu_{t+1}\mbox{ if }k\le t;\\
t-i+\mu_i-\lambda_{t+1}\mbox{ if }k> t.
\end{array}
\right.
$$

\begin{theorem}\label{theorem:commpoly:1}
Let $\lambda\in X^+(n)$, $\mu\in X^+(n-1)$, $\mu\longleftarrow\lambda$,
$1\le i<n$, $1\le d<p$ and $M\subset(i..n)$.
Then $T_{i,n}^{(d)}(M,1)f_{\mu,\lambda}$ is a nonzero $U(n-1)$-high weight vector
if and only if there is \mbox{an injection $\gamma:M\to[i..n)\times[1..d]$
such that}
\begin{enumerate}
\item\label{theorem:commpoly:1:part:1}
      $\gamma_1(m)\ge m$ and $B^{\mu,\lambda}(m,\gamma_1(m))\equiv d-\gamma_2(m)\pmod p$
      for any $m\in M$;
\item\label{theorem:commpoly:1:part:2}
      $B^{\mu,\lambda}(i,t)\not\equiv d-s\pmod p$ for any $(t,s)\in[i..n)\times[1..d]\setminus\im\gamma$,
\end{enumerate}
where $\gamma(m)=(\gamma_1(m),\gamma_2(m))$.
\end{theorem}
\Proof We put $v{:=}T_{i,n}^{(d)}(M,1)f_{\mu,\lambda}$,
$a_t{:=}\sum_{s=1}^t(\lambda_s{-}\mu_s)$,
$A{:=}(a_1,{\ldots},a_{n-1})$ and $\Lambda_t:=[t..n)\times[1..d]$.
By Lemma~\ref{lemma:coeff:1}, we have
$$
\left(\prod_{t=1}^{n-1}\E_t^{(a_t+d\delta_{t\in[i..n)})}\right)
\T_{i,n}^{(d)}(M,1)\equiv \\
\rho^{(A)}(i,n,(n^d),(0^d),M,1)
\left(\prod_{t=1}^{n-1}\E_t^{(a_t)}\right)
$$
modulo $\bar\J^{(A)}$. By Lemmas~\ref{lemma:int:1} and~\ref{lemma:coeff:0.5},
both sides of this equivalence belong to $\U(n)$.
Therefore, applying Lemma~\ref{lemma:algebras:2}, we obtain that this
equivalence holds modulo $\J^{(A)}$.
Tensoring the above equivalence with $1_\K$ gives
$$
\left(\prod_{t=1}^{n-1}E_t^{(a_t+d\delta_{t\in[i..n)})}\right)
T_{i,n}^{(d)}(M,1)\equiv \\
\xi^{(A)}(i,n,(n^d),(0^d),M,1)
\left(\prod_{t=1}^{n-1}E_t^{(a_t)}\right)
$$
modulo $J^{(A)}$.
Since $J^{(A)}f_{\mu,\lambda}=0$, we can multiply this equivalence by $f_{\mu,\lambda}$
on the right and apply Proposition~\ref{proposition:primvect:1}
to the right-hand side. Hence
\begin{equation}\label{equation:commpoly:2}
\c(v)=\pi_\lambda(\xi^{(A)}(i,n,(n^d),(0^d),M,1)).
\end{equation}

{\it``If part''}. Lemma~\ref{lemma:commpoly:1} and
condition~\ref{theorem:commpoly:1:part:1} give
$\pi_\lambda(\xi^{(A)}(i,n,(n^d),(0^d),M,1))$\linebreak$=\prod_{(t,s)\in \Lambda_i\setminus\im\gamma}
(B^{\mu,\lambda}(i,t)-d+s)+p\Z$.
This element is not equal to zero by virtue of
condition~\ref{theorem:commpoly:1:part:2}.
Thus by~(\ref{equation:commpoly:2}), we have proved $\c(v)\ne0$.

Let $V$ be the subspace of $\nabla_n(\lambda)$ spanned by all
vectors of the form $XT_{l,n}^{(d)}(M\cap(l..n),1)f_{\mu,\lambda}$,
where $X\in U^{-,0}(i,l)$
and $l\in M$. We claim that
{\renewcommand{\labelenumi}{{\rm \theenumi}}
\renewcommand{\theenumi}{{\rm(\alph{enumi})}}
\begin{enumerate}
\item\label{theorem:commpoly:1:property:1}
      $\c(u)=0$ for any $u$ of a $U(n)$-weight space of $V$;
\item\label{theorem:commpoly:1:property:2}
      $u\in V$ for any $u$ such that $u\ne v$ and $v\dasharrow u$.
\end{enumerate}}
These properties, once proved, will immediately imply that $v$ is
a nonzero $U(n-1)$-high weight vector by virtue of
Theorem~\ref{theorem:primvect:1}\ref{theorem:primvect:1:part:2}.

To prove~\ref{theorem:commpoly:1:property:1}, we can restrict ourselves
to the case $u=XT_{l,n}^{(d)}(M\cap(l..n),1)f_{\mu,\lambda}$,
where $X$ is an element of $U^{-,0}(i,l)$ having $U(n)$-weight
$\sigma=-m_1\alpha_1-\cdots-m_{n-1}\alpha_{n-1}$
and $l\in M$.
We put $\bar a_t:=a_t+d\delta_{t\in[l..n)}$ and
$\bar A:=(\bar a_1,\ldots,\bar a_{n-1})$.
By Lemma~\ref{lemma:coeff:2}, we get
$
\left(\prod_{t=1}^{n-1}E_t^{(
\bar a_t+m_t
)}\right)X{\equiv }
\xi\left(\prod_{t=1}^{n-1}E_t^{(
\bar a_t
)}\right)
\pmod{J^{(\bar A)}}
$
for some $\xi\in U^0(n)$.
Since $\lambda$ is the highest weight of $\nabla_n(\lambda)$,
we have $J^{(\bar A)}T_{l,n}^{(d)}(M\cap(l..n),1)f_{\mu,\lambda}=0$.
Hence
$$
\displaystyle \left(\prod_{t=1}^{n-1}E_t^{(
\bar a_t+m_t\!
)}\right)\!
XT_{l,n}^{(d)}(M\cap(l..n),1)f_{\mu,\lambda}{=}\xi\!\left(\prod_{t=1}^{n-1}E_t^{(
\bar a_t\!
)}\right)\!
T_{l,n}^{(d)}(M\cap(l..n),1)f_{\mu,\lambda}.
$$

Lemma~\ref{lemma:coeff:1}, yields
$$
\begin{array}{l}
\displaystyle \left(\prod_{t=1}^{n-1}\E_t^{(
\bar a_t
)}\right)
\T_{l,n}^{(d)}(M\cap(l..n),1)\\[12pt]
\displaystyle \equiv \rho^{(A)}(l,n,(n^d),(0^d),M\cap(l..n),1)
\left(\prod_{t=1}^{n-1}\E_t^{(a_t)}\right)
\pmod{\bar\J^{(A)}}.
\end{array}
$$
Again applying Lemmas~\ref{lemma:int:1},~\ref{lemma:coeff:0.5}
and~\ref{lemma:algebras:2}, we obtain that this
equivalence holds modulo $\J^{(A)}$.
Therefore, tensoring this equivalence with $1_\K$,
multiplying both sides by $f_{\mu,\lambda}$ on the right and
applying Proposition~\ref{proposition:primvect:1}, we get
$$
\displaystyle \left(\prod_{t=1}^{n-1}E_t^{(
\bar a_t
)}\right)
T_{l,n}^{(d)}(M\cap(l..n),1)f_{\mu,\lambda}{=}\pi_\lambda(\xi^{(A)}(l,n,(n^d),(0^d),M\cap(l..n),1))f_\lambda.
$$
By Lemma~\ref{lemma:commpoly:1} and
condition~\ref{theorem:commpoly:1:part:1}, the right-hand side of the last equality equals
$\prod_{(t,s)\in \Lambda_l\setminus\gamma(M\cap(l..n))}(B^{\mu,\lambda}(l,t)-d+s)f_\lambda=0$.
Hence $\c(u)=0$.

To prove~\ref{theorem:commpoly:1:property:2},
we shall prove at first that $V$ is closed under
multiplication by elements $E_t^{(r)}$, where
$1\le t<n-1$ and $r\ge0$. We use induction on $r$.
Since the case $r=0$ is obvious, we assume that $r>0$ and
that $V$ is closed under multiplication by $E_t^{(r-1)}$.
It suffices to prove
$E_t^{(r)}u\in V$
for $u=XT_{l,n}^{(d)}(M\cap(l..n),1)f_{\mu,\lambda}$,
where $l\in M$ and $X\in U^{-,0}(i,l)$.
Clearly, this is true if $t\notin[i..n)$ as in that case $E_t^{(r)}u=0$.
In the case $t\in[i..l)$, the result
follows from Lemma~\ref{lemma:commpoly:3}. Finally, let $t\in[l..n)$.
Then we have $E_t^{(r)}u=YE_t^{(r)}T_{l,n}^{(d)}(M\cap(l..n),1)f_{\mu,\lambda}$
for some $Y\in U^{-,0}(i,l)$. It follows from Lemma~\ref{lemma:commpoly:3}
that $E_t^{(r)}T_{l,n}^{(d)}(M\cap(l..n),1)f_{\mu,\lambda}=0$ if $r>d$.
Thus we assume that $r\le d<p$. Taking into account
the formula $E_t^{(r)}=r^{-1}E_t^{(r-1)}E_t$ and the inductive hypothesis,
we can assume that $r=1$.
In that case, the required result follows from Lemma~\ref{lemma:low5:3},
Lemma~\ref{lemma:int:1}, Corollary~\ref{corollary:int:1}
and Lemma~\ref{lemma:algebras:2}.

Now we are going to prove that
$E_t^{(r)}v\in V$ for $1\le t<n-1$ and $r>0$.
If $r>d$ then $E_t^{(r)}v=0$ by Lemma~\ref{lemma:commpoly:3}.
Therefore we assume that $r\le d<p$.
Since $E_t^{(r)}{=}r^{-1}E_t^{(r-1)}E_t$ and $E_t^{(r-1)}V\subset V$,
it remains to show that $E_tv\in V$.
This follows from Lemma~\ref{lemma:low5:3},
Lemma~\ref{lemma:int:1}, Corollary~\ref{corollary:int:1}
and Lemma~\ref{lemma:algebras:2}.

{\it ``Only if part''}.
We shall prove by downward induction on $q{\in}M{\cup}\{n\}$ that
\begin{equation}\label{equation:commpoly:4.5}
\begin{array}{l}
\mbox{for any $m\in M\cap[q..n)$
there are $t_m\in[m..n)$ and $s_m\in[1..d]$}\\
\mbox{such that the pairs $(t_m,s_m)$ are mutually distinct and}\\
\mbox{$B^{\mu,\lambda}(m,t_m)\equiv d-s_m\!\!\!\!\!\pmod p$ for all $m\in M\cap[q..n)$.}
\end{array}
\end{equation}
Suppose that $M\ne\emptyset$, $q\in(M\cup\{n\})\setminus\{\min M\}$ and
for all $m\in M\cap[q..n)$ the numbers $t_m$ and $s_m$
have already been defined. Let $l$ denote the element of $M$
directly preceding $q$. We must define $t_l$ and $s_l$.

By Corollary~\ref{corollary:commpoly:1}, the polynomial
$\rho^{(A)}(i,n,(n^d),(0^{d+1}),M,1)$ belongs to the ideal of
$\mathcal U^0(i,n)$ generated by the polynomials
$$
\prod_{r=0}^k\tfrac1d\rho^{(A)}\left(i_r,i_{r+1},(i_{r+1}-1,i_{r+1}^{d-1}),(0^d,d),M\cap(i_r..i_{r+1}),\tfrac1{\zeta-d}\right),
$$
where
$
i=i_0<i_1<\cdots<i_k<i_{k+1}=l$, $k\ge 0$ and $i_1,\ldots,i_k\in M$.
By~(\ref{equation:commpoly:2}) and
Theorem~\ref{theorem:primvect:1}\ref{theorem:primvect:1:part:2},
we get $\pi_\lambda(\xi^{(A)}(i,n,(n^d),(0^d),M,1)){\ne}0$.
Therefore, there exist integers $i_0,i_1,\ldots,i_k,i_{k+1}$
satisfying the above conditions such that
\begin{equation}\label{equation:commpoly:4.875}
\prod_{r=0}^k\pi_\lambda\left(\xi^{(A)}\!\!\left(i_r,i_{r+1},(i_{r+1}{-}1,i_{r+1}^{d-1}),(0^d,d),M\cap(i_r..i_{r+1}),\tfrac1{\zeta-d}\right)\right){\ne}0.
\end{equation}
We put
$K^{(s)}{:=}(i_s-1,i_s^{d-1})$ for
$1\le s\le k+1$,
$L^{(s)}:=(0^d,\delta_{i_s=i_{s+1}-1})$ for
$1{\le}s{\le}k$
and
$L^{(k+1)}{:=}(0^d)$.
Lemma~\ref{lemma:low5:3}, Lemma~\ref{lemma:low5:1}\ref{lemma:low5:1:part:2}
and Corollary~\ref{corollary:int:1} yeild
\begin{equation}\label{equation:commpoly:5}
\E_{i_{k+1}-1}\cdots\E_{i_1-1}\T_{i,n}^{(d)}(M,1)\equiv
\X_1\cdots\X_{k+1}\T_{l,n}^{(d)}(M\cap(l..n),1)
\end{equation}
modulo $\bar\I_{\{i_1-1,\ldots,i_{k+1}-1\}}$,
where each $\X_s\in\mathcal U^{-,0}(i_{s-1},i_s)$ and
$$
\X_s\equiv -\E_{i_s}(
K^{(s)},
L^{(s)})\left(\T_{i_{s-1},i_s}^{(d)}\left(M\cap(i_{s-1}..i_s),\tfrac1{\zeta-d}\right)\right)
\pmod{\bar\I_{\{i_s-1\}}}.
$$
Lemma~\ref{lemma:int:1} and Lemma~\ref{lemma:algebras:2} show that
equivalence~(\ref{equation:commpoly:5}) holds modulo\linebreak
$\I_{\{i_1-1,\ldots,i_{k+1}-1\}}$.
Let $X_1,\ldots,X_{k+1}$ denote the images in $U(n)$
of $\X_1,\ldots,$\linebreak$\X_{k+1}$ respectively.
Then it follows from~(\ref{equation:commpoly:5}) that
\begin{equation}\label{equation:commpoly:6}
0=E_{i_{k+1}-1}\cdots E_{i_1-1}v=
X_1\cdots X_{k+1}T_{l,n}^{(d)}(M\cap(l..n),1)f_{\mu,\lambda}.
\end{equation}

Let $a^{(s)}_t:=a_t+d\delta_{t\in[i_s..n)}-\delta_{t\in\{i_{s+1}-1,\ldots,i_{k+1}-1\}}$
and $A^{(s)}:=(a^{(s)}_1,\ldots,a^{(s)}_{n-1})$.
By Lemmas~\ref{lemma:coeff:1},~\ref{lemma:coeff:0.5}
and~\ref{lemma:algebras:2}, for any $s=1,\ldots,k+1$, we get
$$
\begin{array}{l}
\displaystyle \left(\prod_{t=1}^{n-1}\E_t^{(a^{(s-1)}_t)}\right)\X_s\equiv
-\rho^{(A^{(s)})}\left(i_{s-1},i_s,K^{(s)},L^{(s)},M\cap(i_{s-1}..i_s),\tfrac1{\zeta-d}\right)\times\\[12pt]
\displaystyle \times\left(\prod_{t=1}^{n-1}\E_t^{(a_t^{(s)})}\right)
\pmod{\J^{(A^{(s)})}+\I_{\{i_s-1\}}}.
\end{array}
$$
It is elementary to see that
$
(\J^{(A^{(s)})}{+}\I_{\{i_s-1\}})
\X_{s+1}{\cdots}\X_{k+1}\T_{l,n}^{(d)}(M\cap(l..n),1){\subset}
$
\linebreak
$
\J^{(A)}+\I_{\{i_s-1\}}$.
Therefore, applying Lemma~\ref{lemma:coeff:3}, we have
$$
\begin{array}{l}
\displaystyle \left(\prod_{t=1}^{n-1}\E_t^{(a^{(0)}_t)}\right)
\X_1\cdots\X_{k+1}\T_{l,n}^{(d)}(M\cap(l..n),1)\equiv (-1)^{k+1}\times\\[6pt]
\displaystyle \times\left(\prod_{s=1}^{k+1}\rho^{(A)}\left(i_{s-1},i_s,K^{(s)},(0^d,d),M\cap(i_{s-1}..i_s),\tfrac1{\zeta-d}\right)\right)\times\\[6pt]
\times\displaystyle \left(\prod_{t=1}^{n-1}\E_t^{(a_t+d\delta_{t\in[l..n)})}\right)\T_{l,n}^{(d)}(M\cap(l..n),1)
\pmod{\J^{(A)}+\I_{\{i_1-1,\ldots,i_{k+1}-1\}}}.
\end{array}
$$
Tensoring this with $1_\K$, multiplying by $f_{\mu,\lambda}$ and
applying~(\ref{equation:commpoly:6}),~(\ref{equation:commpoly:4.875})
and Lemma~\ref{lemma:commpoly:1} and Proposition~\ref{proposition:primvect:1}, we get
$$
\begin{array}{l}
\displaystyle 0=\left(\prod_{t=1}^{n-1}E_t^{(a_t+d\delta_{t\in[l..n)})}\right)T_{l,n}^{(d)}(M\cap(l..n),1)f_{\mu,\lambda}\\[12pt]
\displaystyle =\pi_\lambda(\xi^{(A)}(l,n,(n^d),(0^d),M\cap(l..n),1))f_\lambda\\[12pt]
\displaystyle =\prod\Bigl\{B^{\mu,\lambda}(l,t)-d+s:
\displaystyle (t,s)\in
\Lambda_l
\setminus
\{(t_m,s_m):m\in M\cap[q..n)\}
\Bigr\}f_\lambda.
\end{array}
$$
This formula allows us to choose $t_l$ and $s_l$ as required.

Applying~(\ref{equation:commpoly:4.5}) for $q=\min M\cup\{n\}$,
we define $\gamma(m):=(t_m,s_m)$ for $m\in M$. For this $\gamma$
condition~\ref{theorem:commpoly:1:part:1} is clearly satisfied.
Lemma~\ref{lemma:commpoly:1} now implies
$$
\begin{array}{l}
\displaystyle 0\ne\pi_\lambda(\xi^{(A)}(i,n,(n^d),(0^d),M,1))
\displaystyle =\prod\nolimits_{(t,s)\in \Lambda_i\setminus\im\gamma}
(B^{\mu,\lambda}(i,t){-}d{+}s
)+p\Z,
\end{array}
$$
which gives condition~\ref{theorem:commpoly:1:part:2}.
\endproof

\begin{theorem}\label{theorem:commpoly:2}
Let $\lambda\in X^+(n)$, $\mu\in X^+(n-1)$, $\mu\longleftarrow\lambda$,
$1\le i<j<n$, $1\le d<p$ and $M\subset(i..j)$.
Then $T_{i,j}^{(d)}(M,1)f_{\mu,\lambda}$ is a nonzero $U(n-1)$-high weight vector
if and only if for any sequence of integers $K=(k_1,\ldots,k_d)$ such that
$i\le k_1\le\cdots\le k_d\le j$, there exists an injection
$\gamma_K:M\to[i..j)\times[1..d]$ such that
\begin{enumerate}
\item\label{theorem:commpoly:2:part:1}
      $\gamma_{K,1}(m)\ge m$ and $B^{\mu,\lambda,k_{\gamma_{K,2}(m)}}(m,\gamma_{K,1}(m))\equiv d-\gamma_{K,2}(m)\pmod p$
      for any $m\in M$;
\item\label{theorem:commpoly:2:part:2}
      if $K\ne(j^d)$ then the exits a pair
      $(t,s)\in[i..j)\times[1..d]\setminus\im\gamma_K$ such that
      $B^{\mu,\lambda,k_s}(i,t)\equiv d-s\pmod p$;
\item\label{theorem:commpoly:2:part:3}
      $B^{\mu,\lambda}(i,t)\not\equiv d-s\pmod p$ for any pair
      $(t,s)\in[i..j)\times[1..d]\setminus\im\gamma_{(j^d)}$,
\end{enumerate}
where $\gamma_K(m)=(\gamma_{K,1}(m),\gamma_{K,2}(m))$.
\end{theorem}
\Proof We put $v{:=}T_{i,j}^{(d)}(M,1)f_{\mu,\lambda}$,
$a_t{:=}\sum_{s=1}^t(\lambda_s{-}\mu_s)$,
$A{:=}(a_1,{\ldots},a_{n-1})$ and $\Lambda_t:=[t..j)\times[1..d]$.

{\it``If part''}. Using condition~\ref{theorem:commpoly:2:part:1}
for $K=(j^d)$ and condition~\ref{theorem:commpoly:2:part:3}, we
prove that $\c(v)\ne0$ similarly to how it was done in
Theorem~\ref{theorem:commpoly:1}.

Let $V$ be the subspace of $\nabla_n(\lambda)$ spanned by
vectors of the form
{\renewcommand{\labelenumi}{{\rm \theenumi}}
\renewcommand{\theenumi}{{\rm(\arabic{enumi})}}
\begin{itemize}
\item\label{vector:1}
      $E_j(K,(0^d))\bigl(T_{i,j}^{(d)}(M,1)\bigl)f_{\mu,\lambda}$,
      where $K$ is a weakly increasing sequence distinct from $(j^d)$
      with entries from $[i..j]$;
\item\label{vector:2}
      $X E_j(K,(0^d))\bigl(T_{l,j}^{(d)}(M\cap(l..j),1)\bigl)f_{\mu,\lambda}$,
      where $K$ is a weakly increasing sequence
      with entries from $[l..j]$, $X\in U^{-,0}(i,l)$ and $l\in M$.
\end{itemize}}

Properties~\ref{theorem:commpoly:1:property:1} and~\ref{theorem:commpoly:1:property:2}
from the ``if-part'' of Theorem~\ref{theorem:commpoly:1} also hold in the present situation.
The only alteration that must be made in the proofs is that we now additionally
apply Lemma~\ref{lemma:low5:4}.

{\it ``Only if part''}.
We fix a sequence of integers $K=(k_1,{\ldots},k_d)$ such that
$i{\le}$\linebreak $k_1\le{\cdots}{\le}k_d{\le}j$.
We shall prove by downward induction on $q{\in}M{\cup}\{j\}$ that
\begin{equation}\label{equation:commpoly:12}
\begin{array}{l}
\mbox{for any $m\in M\cap[q..j)$ there are $t_m\in[m..j)$ and $s_m\in[1..d]$}\\
\mbox{such that the pairs $(t_m,s_m)$ are mutually distinct and}\\
\mbox{$B^{\mu,\lambda,k_{s_m}}(m,t_m)\equiv d-s_m\!\!\!\!\!\pmod p$ for any $m\in M\cap[q..j)$.}
\end{array}
\end{equation}
Suppose that $M\ne\emptyset$,
$q\in(M\cup\{j\})\setminus\{\min M\}$
and for all
$m\in M\cap[q..j)$
the numbers $t_m$ and $s_m$
have already been defined. Let $l$ denote the element of $M$
directly preceding $q$. We must define $t_l$ and $s_l$.

Similarly to~(\ref{equation:commpoly:4.875}) we prove that there are integers
$i_1,\ldots,i_k$ belonging to $M$ such that $k\ge 0$, $i=i_0<i_1<\cdots<i_k<i_{k+1}=l$ and
\begin{equation}\label{equation:commpoly:14}
\prod_{r=0}^k\pi_\lambda\!\left(\xi^{(A)}\!\left(i_r,i_{r+1},(i_{r+1}{-}1,i_{r+1}^{d-1}),(0^d,d),M\cap(i_r..i_{r+1}),\tfrac1{\zeta-d}\right)\right){\ne}0.
\end{equation}

We put for brevity $K^{(s)}:=(i_s-1,i_s^{d-1})$ for $s=1,\ldots,k+1$,
$L^{(s)}:=(0^d,\delta_{i_s=i_{s+1}-1})$ for $s=1,\ldots,k$ and
$L^{(k+1)}:=(0^d,b)$, where $b$ is the number of entries of $K$
not greater than $l$.
Lemmas~\ref{lemma:low5:3},~\ref{lemma:low5:1}\ref{lemma:low5:1:part:2}
and Corollary~\ref{corollary:int:1} yield
\begin{equation}\label{equation:commpoly:15}
\begin{array}{l}
\E(k_d^{(l,j)},j)\cdots\E(k_1^{(l,j)},j)
\E_{i_{k+1}-1}\cdots\E_{i_1-1}\T_{i,j}^{(d)}(M,1)\equiv \\
\X_1\cdots\X_{k+1}
\E_j(K^{(l,j)},(0^d))\bigl(\T_{l,j}^{(d)}(M\cap(l..j),1)\bigl)
\end{array}
\end{equation}
modulo $\bar\I_{\{i_1-1,\ldots,i_{k+1}-1\}}$,
where each $\X_s\in\mathcal U^{-,0}(i_{s-1},i_s)$ and
$$
\X_s\equiv -\E_{i_s}(K^{(s)},L^{(s)})
\left(\T_{i_{s-1},i_s}^{(d)}\left(M\cap(i_{s-1}..i_s),\tfrac1{\zeta-d}\right)\right)
\pmod{\bar\I_{\{i_s-1\}}}.
$$
Lemmas~\ref{lemma:int:1} and~\ref{lemma:algebras:2} show
that~(\ref{equation:commpoly:15}) holds modulo
$\I_{\{i_1-1,\ldots,i_{k+1}-1\}}$.
Let $X_1,\ldots,X_{k+1}$ denote the images in $U(n)$ of $\X_1,\ldots,\X_{k+1}$
respectively.

Let $a^{(s)}_t:=a_t+\sum_{h=1}^d\delta_{t\in[i_s..k_h^{(l,j)})}-\delta_{t\in\{i_{s+1}-1,\ldots,i_{k+1}-1\}}$
and $A^{(s)}:=(a^{(s)}_1,$\linebreak$\ldots,a^{(s)}_{n-1})$.
By Lemmas~\ref{lemma:coeff:1},~\ref{lemma:coeff:0.5}
and~\ref{lemma:algebras:2}, for any $s=1,\ldots,k+1$, we get
$$
\begin{array}{l}
\displaystyle \left(\prod_{t=1}^{n-1}\E_t^{(a^{(s-1)}_t)}\right)\X_s\equiv
-\rho^{(A^{(s)})}\left(i_{s-1},i_s,K^{(s)},L^{(s)},M\cap(i_{s-1}..i_s),\tfrac1{\zeta-d}\right)\times\\[12pt]
\displaystyle \times\left(\prod_{t=1}^{n-1}\E_t^{(a_t^{(s)})}\right)
\pmod{\J^{(A^{(s)})}+\I_{\{i_s-1\}}}.
\end{array}
$$

It is elementary to verify that
$
\left(\J^{(A^{(s)})}+\I_{\{i_s-1\}}\right)
\X_{s+1}\cdots\X_{k+1}\times$\linebreak
$
{\times}\E_j(K^{(l,j)},(0^d))\bigl(\T_{l,j}^{(d)}(M\cap(l..j),1)\bigl)
{\subset}\J^{(A)}{+}\I_{[i..j)}.
$
Hence Lemma~\ref{lemma:coeff:3} yields
$$
\begin{array}{l}
\displaystyle \left(\prod_{t=1}^{n-1}\E_t^{(a^{(0)}_t)}\right)
\X_1\cdots\X_{k+1}
\E_j(K^{(l,j)},(0^d))\bigl(
\T_{l,j}^{(d)}(M\cap(l..j),1)
\bigl)
\equiv (-1)^{k+1}\times\\[6pt]
\displaystyle \times\left(\prod_{s=1}^{k+1}\rho^{(A)}\left(i_{s-1},i_s,K^{(s)},(0^d,d),M\cap(i_{s-1}..i_s),\tfrac1{\zeta-d}\right)\right)\times\\[6pt]
\times\displaystyle
\left(\prod_{t=1}^{n-1}\E_t^{(a_t^{(k+1)})}\right)
\E_j(K^{(l,j)},(0^d))\bigl(
\T_{l,j}^{(d)}(M\cap(l..j),1)
\bigl)
\!
\!
\!
\pmod{\J^{(A)}+\I_{[i..j)}}.
\end{array}
$$
By Lemmas~\ref{lemma:coeff:1},~\ref{lemma:coeff:0.5} and~\ref{lemma:algebras:2}
the last product of this equivalence equals
$\rho^{(A)}(l,j,K^{(l,j)},(0^d),M\cap(l..j),1)\prod_{t=1}^{n-1}\E_t^{(a_t)}$
modulo $\pmod{\J^{(A)}+\I_{[i..j)}}$.
We note that $B^{\mu,\lambda,k_s^{(l,j)}}(q,t)=B^{\mu,\lambda,k_s}(q,t)$ for $t\ge l$.
Applying this remark, the above equivalences, equivalence~(\ref{equation:commpoly:15}),
inequality~(\ref{equation:commpoly:14}), Theorem~\ref{theorem:primvect:1}\ref{theorem:primvect:1:part:2}
and Lemma~\ref{lemma:commpoly:1}, we get
$$
\begin{array}{l}
\displaystyle 0=\pi_\lambda(\xi^{(A)}(l,j,K^{(l,j)},(0^d),M\cap(l..j),1))f_\lambda=\\[12pt]
\displaystyle =\prod\Bigl\{B^{\mu,\lambda,k_s}(l,t)-d+s:
\displaystyle (t,s)\in \Lambda_l\setminus
\{(t_m,s_m):m\in M\cap[q..j)\}
\Bigr\}f_\lambda.
\end{array}
$$
This formula allows us to choose $t_l$ and $s_l$ as required.

Applying~(\ref{equation:commpoly:12}) for $q=\min M\cup\{j\}$,
we define $\gamma_K(m):=(t_m,s_m)$ for $m\in M$. For this $\gamma_K$
condition~\ref{theorem:commpoly:2:part:1} is clearly satisfied.
Lemma~\ref{lemma:commpoly:1} now implies
$$
\begin{array}{l}
\displaystyle\left(\prod_{t=1}^{n-1}E_t^{(\bar a_t)}\right)
(E_{k_d}\cdots E_{j-1})\cdots(E_{k_1}\cdots E_{j-1})
T_{i,j}^{(d)}(M,1)f_{\mu,\lambda}=\\[6pt]
\displaystyle\pi_\lambda(\xi^{(A)}(i,j,K,(0^d),M,1))f_\lambda
{=}\!\left(\prod\nolimits_{(t,s)\in \Lambda_i\setminus\im\gamma_K}
(B^{\mu,\lambda,k_s}(i,t)-d+s)\right)\!f_\lambda,
\end{array}
$$
where $\bar a_t=a_t+\sum_{s=1}^d\delta_{t\in[i..k_s)}$.
By Theorem~\ref{theorem:primvect:1}\ref{theorem:primvect:1:part:2}, the left-hand side of
the above formula equals zero if and only if $K\ne(j^d)$.
Hence conditions~\ref{theorem:commpoly:2:part:2} and~\ref{theorem:commpoly:2:part:3}
easily follow. \endproof

Simpler criterions will appear if we consider the problem of
existence of $M$ such that $T_{i,j}^{(d)}(M,1)f_{\mu,\lambda}$
is a nonzero $U(n-1)$-high weight vector in Theorem~\ref{theorem:commpoly:1}
or in Theorem~\ref{theorem:commpoly:2} in the cases $j=n$ or $j<n$, respectively.

For what follows, let us recall the definition of the sets
$\mathfrak C^\mu(i,j)$, $\mathfrak X^\mu_d(i,j)$, $\mathfrak X^{\mu,\lambda}_d(i,j)$
and Definition~\ref{def:introduction:1} (see the introduction).
We also use the notation $C^\mu(i,t)=t-i+\mu_i-\mu_t$.\label{Cmu}

\begin{theorem}\label{theorem:crit:0.5}
Let $1\le i<n$, $1\le d<p$, $\mu\in X^+(n-1)$,
$\lambda\in X^+(n)$ and $\mu\longleftarrow\lambda$.
There exists $M\subset(i..n)$ such that $T_{i,n}^{(d)}(M,1)f_{\mu,\lambda}$ is a nonzero
$U(n-1)$-high weight vector if and only if there exists an injection
$\epsilon:\mathfrak X^{\mu,\lambda}_d(i,n)\to\mathfrak C^\mu(i,n)$ weakly decreasing
w.r.t. the first coordinate.
\end{theorem}
\Proof
At first, we suppose that such $M$ exists.
Let $\gamma:M\to[i..n)\times[1..d]$ be
an injection satisfying conditions~\ref{theorem:commpoly:1:part:1}
and~\ref{theorem:commpoly:1:part:2} of Theorem~\ref{theorem:commpoly:1}.
Condition~\ref{theorem:commpoly:1:part:2} immediately implies
$\mathfrak X^{\mu,\lambda}_d(i,n)\subset\im \gamma$.
Now let $\gamma(m)\in\mathfrak X^{\mu,\lambda}_d(i,n)$ for some $m\in M$.
Applying condition~\ref{theorem:commpoly:1:part:1}, we easily get
$C^\mu(i,m)\equiv 0\pmod p$ and $m\in\mathfrak C^\mu(i,n)$.
Now we can take for $\epsilon$ the map with domain $\mathfrak X^{\mu,\lambda}_d(i,n)$
partially inverting $\gamma$.

Now let $\epsilon:\mathfrak X^{\mu,\lambda}_d(i,n)\to\mathfrak C^\mu(i,n)$
be an injection as in the formulation of the theorem.
We put $M:=\im\epsilon$ and denote by $\gamma$ the inverse map of $\epsilon$.
Condition~\ref{theorem:commpoly:1:part:2} of Theorem~\ref{theorem:commpoly:1}
follows from $\im\gamma=\mathfrak X^{\mu,\lambda}_d(i,n)$.
Let $m\in M$.
We have $m=\epsilon(t,s)$ for some $(t,s)\in\mathfrak X^{\mu,\lambda}_d(i,n)$
and $C^\mu(i,m)\equiv 0\pmod p$. Since $B^{\mu,\lambda}(i,t)\equiv d-s\pmod p$, we have
$B^{\mu,\lambda}(m,t)\equiv d-s\pmod p$, which can be reformulated as
$B^{\mu,\lambda}(m,\gamma_1(m))\equiv d-\gamma_2(m)\pmod p$ to conform to
condition~\ref{theorem:commpoly:1:part:1} of Theorem~\ref{theorem:commpoly:1}.
\endproof

Before proceeding further, we make a remark on partially ordered sets.
Let $X$ be a finite set with nonstrict partial order $\preccurlyeq$.
We put $\cone(x):=\{y\in X:x\preccurlyeq y\}$\label{conex}
for any $x\in X$
and $\cone(S):=\bigcup_{x\in S}\cone(x)$\label{coneS}
for any $S\subset X$.
A map $\alpha:A\to B$, where $A,B\subset X$, is called
{\it weakly increasing} if $x\preccurlyeq\alpha(x)$ for any $x\in A$.

\begin{proposition}\label{proposition:crit:1}
Let $A,B\subset X$. There exists a weakly increasing injection from
$A$ to $B$ if and only if
$|\cone(S)\cap A|\le |\cone(S)\cap B|$ for any $S\subset X$.
\end{proposition}

For the rest of the paper, we fix integers $i,j,d$ and weights $\mu\in X^+(n-1)$,
$\lambda\in X^+(n)$ such that $1\le i<j<n$, $1\le d<p$ and $\mu\longleftarrow\lambda$.
We put $X:=[i..j)\times[1..d]$. For a sequence of integers
$K=(k_1,\ldots,k_d)$ such that $i\le k_1\le\cdots\le k_d\le j$, we put
$
\mathfrak X^{\mu,\lambda,K}_d(i,j):={\{} (t,s)\in X:B^{\mu,\lambda,k_s}(i,t)\equiv d-s\pmod p {\}}.
$\label{XmulmK}
Using $K$, we define the following subsets of $X$:
$Y_K:=\{(t,s)\in X:t<k_s\}$\label{YK}
and $Z_K:=\{(t,s)\in X:t\ge k_s\}$.\label{ZK}
Clearly $X=Y_K\sqcup Z_K$.
We have
\begin{equation}\label{equation:crit:1}
\mathfrak X^{\mu,\lambda,K}_d(i,j)=
\Bigl(\mathfrak X^{\mu,\lambda}_d(i,j)\cap Y_K\Bigr)\sqcup
\Bigl(\mathfrak X^\mu_d(i,j)\cap Z_K\Bigr).
\end{equation}

\begin{theorem}\label{theorem:crit:1}
Let $1\le i<j<n$, $1\le d<p$, $\mu\in X^+(n-1)$,
$\lambda\in X^+(n)$ and $\mu\longleftarrow\lambda$.

\begin{enumerate}
\item\label{theorem:crit:1:part:1}
There exists $M\subset(i..j)$ such that $T_{i,j}^{(d)}(M,1)f_{\mu,\lambda}$
is a nonzero $U(n-1)$-high weight vector if and only if
there exists an injection $\epsilon:\mathfrak X^{\mu,\lambda}_d(i,j)\to\mathfrak C^\mu(i,j)$
weakly decreasing w.r.t. the first coordinate and
for any sequence of integers $K=(k_1,\ldots,k_d)$ such that
$i\le k_1\le\cdots\le k_d\le j$ and $K\ne(j^d)$ there exists
an injection $\theta_K:\{i\}\cup\im\epsilon\to\mathfrak X^{\mu,\lambda,K}_d(i,j)$
weakly increasing w.r.t. the first coordinate.

\item\label{theorem:crit:1:part:2}
There exists $M\subset(i..j)$ such that
$T_{i,j}^{(d)}(M,1)f_{\mu,\lambda}$ is a nonzero $U(n-1)$-high weight vector
if and only if $(j-1,1)\in\mathfrak X^\mu_d(i,j)$,
$(j-1,1)\notin\mathfrak X^{\mu,\lambda}_d(i,j)$,
there exists an injection
$\epsilon:\mathfrak X^{\mu,\lambda}_d(i,j)\to\mathfrak C^\mu(i,j)$
weakly decreasing w.r.t. the first coordinate
and an injection $\tau:\mathfrak X^{\mu,\lambda}_d(i,j)\to \mathfrak X^\mu_d(i,j)\setminus\{(j-1,1)\}$
weakly increasing w.r.t. the first coordinate and
weakly decreasing w.r.t. the second coordinate.
\end{enumerate}
\end{theorem}
\Proof
\ref{theorem:crit:1:part:1}
Suppose that the required set $M$ exists and let $\gamma_K:M\to X$ be
injections as in Theorem~\ref{theorem:commpoly:2}.
It follows from condition~\ref{theorem:commpoly:2:part:3} of
Theorem~\ref{theorem:commpoly:2} that
$\mathfrak X^{\mu,\lambda}_d(i,j)\subset\im\gamma_{(j^d)}$.
Now let $\gamma_{(j^d)}(m)\in \mathfrak X^{\mu,\lambda}_d(i,j)$ for some $m\in M$.
Condition~\ref{theorem:commpoly:2:part:1} of Theorem~\ref{theorem:commpoly:2} and
the definition of $\mathfrak X^{\mu,\lambda}_d(i,j)$ imply
$m\in\mathfrak C^\mu(i,j)$.
It is clear now that $\epsilon$ can be chosen so that
$\gamma_{(j^d)}\circ\epsilon=id_{\mathfrak X^{\mu,\lambda}_d(i,j)}$.

For a sequence of integers $K=(k_1,\ldots,k_d)$ such that
$i\le k_1\le\cdots\le k_d\le j$ and $K\ne(j^d)$,
we put $\theta_K(m)=\gamma_K(m)$ for $m\in\im\epsilon$ and
put $\theta_K(i)$ equal to the pair $(t,s)$ that is mentioned in
condition~\ref{theorem:commpoly:2:part:2} of Theorem~\ref{theorem:commpoly:2}.

Let $m\in\im\epsilon$.
Applying the congruence $C^\mu(i,m)\equiv 0\pmod p$ and
condition~\ref{theorem:commpoly:2:part:1} of Theorem~\ref{theorem:commpoly:2},
we get $\theta_K(m)=\gamma_K(m)\in\mathfrak X^{\mu,\lambda,K}_d(i,j)$.

Now suppose, on the contrary, that the required maps
$\epsilon$ and $\theta_K$ exist.
We put $M:=\im\epsilon$, $\gamma_{(j^d)}:=\epsilon^{-1}$ and
$\gamma_K:={\theta_K|}_M$ for $K\ne(j^d)$.
We claim that $M$ and $\gamma_K$ so defined satisfy conditions
\ref{theorem:commpoly:2:part:1}--\ref{theorem:commpoly:2:part:3}
of Theorem~\ref{theorem:commpoly:2}.

Let $m\in M$ and $K=(k_1,\ldots,k_d)$ be a sequence
of integers such that $i\le k_1\le\cdots\le k_d\le j$.
Since $\gamma_K(m)\in\mathfrak X^{\mu,\lambda,K}_d(i,j)$ and
$C^\mu(i,m)=0\pmod p$, we have
$B^{\mu,\lambda,k_{\gamma_{K,2}(m)}}(m,\gamma_{K,1}(m))\equiv d-\gamma_{K,2}(m)\pmod p$,
wich gives condition~\ref{theorem:commpoly:2:part:1}.
To see that condition~\ref{theorem:commpoly:2:part:2} holds,
it suffices to $\theta_K(i)$ for $(t,s)$.
Finally condition~\ref{theorem:commpoly:2:part:3} holds,
since $\im\gamma_{(j^d)}=\mathfrak X^{\mu,\lambda}_d(i,j)$.

\ref{theorem:crit:1:part:2}
Let $\epsilon$ and $\tau$ be injections as in the formulation of the theorem.
For a sequence of integers $K=(k_1,\ldots,k_d)$ such that
$i\le k_1\le\cdots\le k_d\le j$ and $K\ne(j^d)$,
we define $\theta_K:\{i\}\cup\im\epsilon\to\mathfrak X^{\mu,\lambda,K}_d(i,j)$
by
$$
\theta_K(x):=\left\{
\begin{array}{ll}
(j-1,1)&\mbox{ if }x=i;\\
\epsilon^{-1}(x)&\mbox{ if }i<x\mbox{ and }\epsilon^{-1}(x)\in Y_K;\\
\tau(\epsilon^{-1}(x))&\mbox{ if }i<x\mbox{ and }\epsilon^{-1}(x)\in Z_K.
\end{array}
\right.
$$
Clearly, these maps $\epsilon$ and $\theta_K$ satisfy
conditions of part~\ref{theorem:crit:1:part:1} of the current theorem
(see~(\ref{equation:crit:1})).

Now, on the contrary, let $\epsilon$ and $\theta_K$ be injections as
in part~\ref{theorem:crit:1:part:1}.
For any sequence of integers $K=(k_1,\ldots,k_d)$ such that
$i\le k_1\le\cdots\le k_d\le j$ and $K\ne(j^d)$, we have
\begin{equation}\label{equation:crit:2}
|\mathfrak X^{\mu,\lambda,K}_d(i,j)|\ge|\im\theta_K|=|\{i\}\cup\im\epsilon|=
1+|\mathfrak X^{\mu,\lambda}_d(i,j)|.
\end{equation}
Putting $K=(j-1,j^{d-1})$ and applying~(\ref{equation:crit:1}), we get
$$
\begin{array}{l}
|\mathfrak X^{\mu,\lambda}_d(i,j)\cap Y_{(j-1,j^{d-1})}|+
|\mathfrak X^\mu_d(i,j)\cap Z_{(j-1,j^{d-1})}|=
|\mathfrak X^{\mu,\lambda,(j-1,j^{d-1})}_d(i,j)|\\[6pt]
\ge 1+|\mathfrak X^{\mu,\lambda}_d(i,j)\cap Y_{(j-1,j^{d-1})}|+
|\mathfrak X^{\mu,\lambda}_d(i,j)\cap Z_{(j-1,j^{d-1})}|.
\end{array}
$$
Hence
$
|\mathfrak X^\mu_d(i,j)\cap Z_{(j-1,j^{d-1})}|\ge
1+|\mathfrak X^{\mu,\lambda}_d(i,j)\cap Z_{(j-1,j^{d-1})}|.
$
Since $Z_{(j-1,j^{d-1})}=\{(j-1,1)\}$, we have
$(j-1,1)\in\mathfrak X^\mu_d(i,j)$ and
$(j-1,1)\notin\mathfrak X^{\mu,\lambda}_d(i,j)$.

To prove the existence of $\tau$, we apply
Proposition~\ref{proposition:crit:1}.
We intoduce the partial order $\preccurlyeq$ on $X$ by
 $(a,b)\preccurlyeq(x,y)\Leftrightarrow a\le x\mathbin{\&}b\ge y$.
Take a nonempty subset $S\subset X$.
Then we have $\cone(S)=Z_K$ for $K=(k_1,\ldots,k_d)$,
where $k_s=\min\bigl(\cone(S)\cap(\Z\times\{s\})\bigr)\cup\{j\}$.
Notice that $K\ne(j^d)$ and $(j-1,1)\in Z_K$, since $S\ne\emptyset$.
By~(\ref{equation:crit:1}) and~(\ref{equation:crit:2}), we have
$$
\begin{array}{l}
1+|\mathfrak X^{\mu,\lambda}_d(i,j)\cap Y_K|+
  |\mathfrak X^{\mu,\lambda}_d(i,j)\cap Z_K|=
1+|\mathfrak X^{\mu,\lambda}_d(i,j)|\\[6pt]
\le|\mathfrak X^{\mu,\lambda,K}_d(i,j)|=
  |\mathfrak X^{\mu,\lambda}_d(i,j)\cap Y_K|+
  |\mathfrak X^\mu_d(i,j)\cap Z_K|\\[6pt]
  =|\mathfrak X^{\mu,\lambda}_d(i,j)\cap Y_K|+
  |(\mathfrak X^\mu_d(i,j)\setminus{\{}(j-1,1){\}})\cap Z_K|+1.
\end{array}
$$
Hence
$$
|\mathfrak X^{\mu,\lambda}_d(i,j)\cap Z_K|\le
|(\mathfrak X^\mu_d(i,j)\setminus{\{}(j-1,1){\}})\cap Z_K|,
$$
which by Proposition~\ref{proposition:crit:1} yields
the required map $\tau$.
\endproof

\section{Appendix: List of Notations}

\tabcolsep=0pt
\begin{tabular}{p{3.6cm}p{8.5cm}}
$\Z$, $\Q$, $\K$        & sets of integers and rationals, algebraically
                          closed field of characteristic $p>0$;\\[2pt]

$L_n(\lambda)$              & irreducible $\GL_n$-module
                          with highest weight $\lambda$, p.~\pageref{Ln(lm)};\\[2pt]

$X(n)$                  & $\Z^n$, weights, p.~\pageref{X(n)};\\[2pt]

$X^+(n)$                & $\{\lambda\in\Z^n:\lambda_1\ge\cdots\ge\lambda_n\}$, dominant weights, p.~\pageref{domX(n)};\\[2pt]

$\mathfrak C^\mu(i,j)$  & $\{t\in(i..j):t-i+\mu_i-\mu_t\equiv0\!\!\pmod p\}$, p.~\pageref{CXX};\\[2pt]

$\mathfrak X^\mu_d(i,j)$& $\{(t,s)\in [i..j)\times[1..d]:t-i+\mu_i-\mu_{t+1}\equiv d-s\!\!\pmod p\}$, p.~\pageref{CXX};\\[2pt]

$\mathfrak X^{\mu,\lambda}_d(i,j)$ & $\{(t,s)\in [i..j)\times[1..d]:t-i+\mu_i-\lambda_{t+1}\equiv d{-}s\!\!\pmod p \}$, p.~\pageref{CXX};\\[2pt]

$[i..j]$, etc.  &  $\{x\in\Z:i\le x\le j\}$, etc., p.~\pageref{intervals};\\[2pt]

$\epsilon_i$ & $(0,\ldots,0,1,0,\ldots,0)$ having $1$ at position $i$, p.~\pageref{weights};\\[2pt]

$\alpha(i,j)$,\, $\alpha_i$ & $\epsilon_i-\epsilon_j$,\, $\alpha(i,i+1)$ resp., p.~\pageref{weights};\\[2pt]

$x^{\underline n}$& $x\cdots (x-n+1)$ if $n\ge0$ and $1/(x+1)\cdots (x-n)$ if $n<0$, p.~\pageref{descending_factorial};\\[2pt]

$\delta_{\mathcal P}$ & $1$ or $0$ if $\mathcal P$ is true or false respectively, p.~\pageref{delta};\\[2pt]

$U_\Q(n)$ & universal enveloping algebra of ${\mathfrak gl}_\Q(n)$, p.~\pageref{universal_enveloping_algebra};\\[2pt]

$\U(n)$   & hyperalgebra over $\Z$, p.~\pageref{hyperalgebra};\\[2pt]

$UT(n)$   & set of integer strictly upper triangular $n\times n$-matrices
            with nonnegative entries, p.~\pageref{UT(n)};\\[2pt]

$\F^{(N)}$,
$\E^{(N)}$, & $\prod_{1\le a<b\le n}\F_{a,b}^{(N_{a,b})}$,
              $\prod_{1\le a<b\le n} \E_{a,b}^{(N_{a,b})}$ resp., p.~\pageref{FNEN};\\[2pt]

$U^0_\Q(n)$ & $\Q$-subalgebra of $U_\Q(n)$ generated by $\H_1,\ldots,\H_n$, p.~\pageref{U0};\\[2pt]

$\tau_N$    & automorphism of $U^0_\Q(n)$ satisfying~(\ref{equation:algebras:1}), p.~\pageref{tauN};\\[2pt]

$\bar\U(n)$ & right ring of quotients of $U_\Q(n)$ with respect
              to $U^0_\Q(n)\setminus\{0\}$, p.~\pageref{barU(n)};\\[2pt]

$\I_S$,\, $\bar\I_S$ & left ideals generated by $\E^{(r)}_i$ for $r\ge1$,
                       $i\in S$ in $\U(n)$ and $\bar\U(n)$ resp., p.~\pageref{idealsI};\\[2pt]

$\J^{(C)}$,\, $\bar\J^{(C)}$ & left ideals generated by elements of weight
                               with $\alpha_i$-coefficient $>c_i$ for some $i$
                               in $\U(n)$ and $\bar\U(n)$ resp., p.~\pageref{idealsJ};\\[2pt]

$\C(i,j)$,\! $\B(i,j)$ & $j-i+\H_i-\H_j$ and $j-i+\H_i-\H_{j+1}$ resp. p.~\pageref{CB};\\[2pt]

$\S_{i,j}$ & Carter--Lusztig lowering operator, p.~\pageref{S};\\[2pt]

$\theta_i$ & p.~\pageref{theta};\\[2pt]

$\E_j(K,L)$& Definition~\ref{definition:low5:1}, p.~\pageref{definition:low5:1};\\[2pt]

$K^{(i,j)}$, $K^{\{j\}}$ & p.~\pageref{minmax};\\[2pt]

$\T_{i,j}^{(d)}(M,R)$ & Definition~\ref{definition:low5:2}, p.~\pageref{definition:low5:2};\\[2pt]

$\B^{C,k}(i,t)$ & $\B(i,t)+c_{i-1}-c_i+\delta_{t\ge k}(c_{t+1}-c_t)$, p.~\pageref{BCk};\\[2pt]

$\rho^{(C)}(i,j,K,L,M,R)$ & equations~(\ref{equation:coeff:2.5}) and~(\ref{equation:coeff:2.75}), p.~\pageref{equation:coeff:2.75};\\[2pt]

$f_{i,j}^{(d)}(M)$,
$g_{i,j}^{(d)}(M)$ & p.~\pageref{fijgij};\\[2pt]

$G_{i,l}^{(d)}(M,N)$ & p.~\pageref{Gil};\\[2pt]

$U(n)$ & $\U(n)\otimes_\Z\K$, hyperalgebra over $\K$, p.~\pageref{U(n)};\\[2pt]

\end{tabular}

\begin{tabular}{p{3.6cm}p{8.5cm}}

$\pi_\lambda$ & $\K$-algebra homomorphism from $U^0(n)$ to $\K$ taking
            $\binom{H_i}{r}$ to $\binom{\lambda_i}{r}+p\Z$, p.~\pageref{pilm};\\[2pt]

$\nabla_n(\lambda)$ & co-Weyl module with highest weight $\lambda$, p.~\pageref{CoWeyl};\\[2pt]

$f_{\mu,\lambda}$ & $U(n{-}1)$-high weight vector of $\nabla_n(\lambda)$
                of weight $\mu$, p.~\pageref{fmulm};\\[2pt]

$f_\lambda$ & $f_{(\lambda_1,\ldots,\lambda_{n-1}),\lambda}$, p.~\pageref{flm};\\[2pt]

$\c(v)$ & element of $\K$ such that
          $E_1^{(a_1)}\cdots E_{n-1}^{(a_{n-1})}v=\c(v)f_\lambda$
          for $v\in\nabla_n(\lambda)$ of weight $\lambda{-}a_1\alpha_1{-}{\cdots}{-}a_{n-1}\alpha_{n-1}$,
          p.~\pageref{cfv};\\[2pt]

$S_{i,j}$ & image of $\S_{i,j}$ in $U(n)$, p.~\pageref{SinU};\\[2pt]

$T_{i,j}^{(d)}(M,R)$ & image of $\T_{i,j}^{(d)}(M,R)$ in $U(n)$, p.~\pageref{TinU};\\[2pt]

$\xi^{(C)}(i,j,K,L,M,R)$ & image of $\rho^{(C)}(i,j,K,L,M,R)$ in $U(n)$, p.~\pageref{xi};\\[2pt]

$I_S$, $J^{(C)}$ & images of $\I_S$ and $\J^{(C)}$ in $U(n)$ resp., p.~\pageref{IJinU};\\[2pt]

$B^{\mu,\lambda}(i,t)$ & $t-i+\mu_i-\lambda_{t+1}$, p.~\pageref{Bmulm};\\[2pt]

$B^{\mu,\lambda,k}(i,t)$ & $t{-}i{+}\mu_i{-}\mu_{t+1}$ if $k{\le}t$ or $t{-}i{+}\mu_i{-}\lambda_{t+1}$ if $k{>}t$, p.~\pageref{Bmulm};\\[2pt]

$C^\mu(i,t)$ &  $t-i+\mu_i-\mu_t$, p.~\pageref{Cmu};\\[2pt]

$\cone(x)$ & $\{y\in X:x\preccurlyeq y\}$, p.~\pageref{conex};\\[2pt]

$\cone(S)$ & $\bigcup_{x\in S}\cone(x)$, p.~\pageref{coneS};\\[2pt]

$\mathfrak X^{\mu,\lambda,K}_d(i,j)$ & $\{(t,s)\in[i..j)\times[1..d]:B^{\mu,\lambda,k_s}(i,t)\equiv d-s\pmod p\}$, p.~\pageref{XmulmK};\\[2pt]

$Y_K$ &  $\{(t,s)\in[i..j)\times[1..d]:t<k_s\}$, p.~\pageref{YK};\\[2pt]

$Z_K$ & $\{(t,s)\in[i..j)\times[1..d]:t\ge k_s\}$, p.~\pageref{ZK}.

\end{tabular}

\bigskip

{\bf Remark.} See the table on page~\pageref{miscellaneous_rings} for the definitions
of $\U^+(a,b)$, $\U^0(a,b)$, $\U^-(a,b)$, $\U^{-,0}(a,b)$,
$\bar\U^0(a,b)$, $\bar\U^{-,0}(a,b)$,
$\mathcal U^0(a,b)$, $\mathcal U^{-,0}(a,b)$.

\providecommand{\bysame}{\leavevmode\hbox to3em{\hrulefill}\thinspace}
\providecommand{\MR}{\relax\ifhmode\unskip\space\fi MR }
\providecommand{\MRhref}[2]{
  \href{http://www.ams.org/mathscinet-getitem?mr=#1}{#2}
}
\providecommand{\href}[2]{#2}


\begin{thebibliography}{BKS}

\bibitem[CL]{Carter_Lusztig}
R.W.Carter and G.Lusztig, On the modular representations of the general
  linear and symmetric groups, {\em Math. Z.}, {\bf 136} (1974), 193--242.

\bibitem[J]{Jantzen1}
J.C.Jantzen, Representations of algebraic groups, Pure and Applied
  Mathematics, 131, {\em Academic Press, Inc.}, Boston, MA, 1987.

\bibitem[K]{Kleshchev2}
A.S.Kleshchev, Branching rules for modular representations of symmetric
  groups. II, {\em J. Reine Angew. Math.}, {\bf 459} (1995), 163--212.

\bibitem[B]{Brundan_operators}
J.Brundan, Lowering operators for {\rm GL}(n) and quantum {\rm GL}(n), {\em
  Proc. Simposia in Pure Math.}, {\bf 63} (1998), 95--114.

\bibitem[BKS]{Kleshchev_gjs11}
J.Brundan, A.Kleshchev and I.Suprunenko, Semisimple restrictions from ${\rm
  GL}(n)$ to ${\rm GL}(n-1)$, {\em J. reine angew. Math.}, {\bf 500} (1998),
  83--112.


\bibitem[BK]{Kleshchev_bou}
J.Brundan and A.S.Kleshchev, Some remarks on branching rules and tensor
  products for algebraic groups, {\em J. Algebra}, {\bf 217} (1999), n. 1,
  335--351.

\bibitem[S]{Shchigolev14}
V.V.Shchigolev, Iterating lowering operators, {\em J. Pure Appl. Alg., Special
  issue dedicated to Eric M. Friedlander on the occasion of his 60th birthday},
  {\bf 206} (2006), 111--122.

\end{thebibliography}
\end{document}